\documentclass[letterpaper,10pt]{article}

\usepackage{amsfonts,amsmath,verbatim,graphicx, subfigure,longtable}
\usepackage{amssymb,amsthm,amscd}
\usepackage{pstricks}
\usepackage{epstopdf}
\usepackage{color}
\usepackage{transparent}
\usepackage{tikz}
\usetikzlibrary{fit, positioning}

\newtheorem{thm}{Theorem}[section]

\newtheorem{lem}[thm]{Lemma}
\newtheorem{definition}[thm]{Definition}
\newtheorem{assumption}[thm]{Assumption}

\newtheorem{remark}[thm]{Remark}

\newenvironment{pf}{\noindent\emph{Proof \,}}{\mbox{}\qed}

\newcommand{\toL}{\,{\buildrel \mathcal{L} \over \longrightarrow}\,}
\newcommand{\equalL}{\,{\buildrel \mathcal{L} \over =}\,}

\def\wh{\widehat}

\numberwithin{equation}{section}
\def\R{{\mathbb R}}     
\def\P{{\mathbb P}}     
\def\E{{\mathbb E}}     
\def\D{{\mathcal E}}    
\def\F{{\mathcal F}}    
\def\H{{\mathcal H}}    %

\def\X{{\mathfrak X}}   
\def\Y{{\mathcal Y}}   

\def\A{{\mathcal A}}   

\def\eps{\varepsilon}
\def\<{{\langle}}
\def\>{{\rangle}}
\def\1{{\bf 1}}         

\renewcommand{\bar}{\overline}

\begin{document}
\allowdisplaybreaks

\title{\Large \bf
Functional central limit theorem for Brownian particles in domains with Robin boundary condition
\thanks{Research partially supported by NSF Grants DMS-1206276 and DMR-1035196.}
 }

\author{{\bf Zhen-Qing Chen} \quad and \quad {\bf Wai-Tong (Louis) Fan}}
\date{\today}
\maketitle

\begin{abstract}
We rigorously derive non-equilibrium space-time fluctuation for the particle density of a system of reflected diffusions in bounded Lipschitz domains in $\R^d$. The particles are independent and are killed by a time-dependent potential which is asymptotically proportional to the boundary local time. We generalize the functional analytic framework introduced by Kotelenez \cite{pK86, pK88} to deal with time-dependent perturbations. Our proof relies on Dirichlet form method rather than the machineries derived from Kotelenez's sub-martingale inequality. Our result holds for any symmetric reflected diffusion, for any bounded Lipschitz domain and for any dimension $d\geq 1$.
\end{abstract}



\bigskip

\tableofcontents

\section{Introduction}

The goal of this paper is to develop a machinery to overcome some difficulties that
 arise in the study of fluctuations for systems of reflected diffusions (such as reflected Brownian motions) with a singular type of time-dependent killing potential. The primary examples are the systems of annihilating diffusions introduced in \cite{zqCwtF13a} and \cite{zqCwtF14b}, which can be used to model the transport of positive and negative charges in solar cells or  the population dynamics of two segregated species under competition.
 The model in \cite{zqCwtF14b} consists of two families of reflected diffusions confined in two adjacent domains, say two adjacent rectangles $(0,2)\times (0,1)$ and $(0,2)\times (-1,0)$, respectively. These two families of particles (positive and negative charges respectively) annihilate each other at a certain rate when they come close to each other near the interface $(0,2)\times \{0\}$. This interaction models the annihilation, trapping, recombination and separation phenomena of the charges. From the viewpoint of the positive charges, they are themselves reflected diffusions in $(0,2)\times (0,1)$ subject to killing by a time-dependent random potential.

In this paper, we focus our attention to a  one-type particle model which consists of i.i.d. reflected diffusions killed by a deterministic time-dependent potential near the boundary. The following assumption on reflected diffusions is in force throughout this paper:

\begin{assumption}\label{A:GeometricSetting}
Suppose $D\subset \R^d$ is a bounded Lipschitz domain, $\rho\in W^{1,2}(D)\cap C(\bar{D})$ is a strictly positive function, $\textbf{a}=(a^{ij})$ is a symmetric, bounded, uniformly elliptic $d\times d$ matrix-valued function such that $a^{ij}\in W^{1,2}(D)$ for each $i,\,j$. Here $C(\bar{D})$ denotes the space of continuous functions on $\bar{D}$ and $W^{1,2}(D):= \{f\in L^2(D):\;|\nabla f|\in L^2(D)\}$ denotes the usual Sobolev space of order $(1, 2)$.
\end{assumption}

Under Assumption \ref{A:GeometricSetting}, it is well known (see \cite{BH91, zqC93}) that the bilinear form $(\D,\,W^{1,2}(D))$ defined by
\begin{equation}\label{E:DirichletForm_ReflectedDiffusion}
\D(f,g)  = \frac{1}{2}\int_{D} \textbf{a} (x) \nabla f(x) \cdot \nabla g(x)\,\rho(x)\,dx  =  \frac{1}{2}\int_D\sum_{i,j=1}^da^{ij}(x)\,\frac{\partial f}{\partial x_i}(x)\,\frac{\partial g}{\partial x_i}(x)\,\rho(x)\,dx
\end{equation}
is a regular Dirichlet form in $L^2(D,\,\rho(x)dx)$ and hence has an associated Hunt process $X$ (unique in distribution). Furthermore, $X$ is a continuous strong Markov process with symmetrizing measure $\rho$ and has infinitesimal generator $\A  = \dfrac{1}{2\,\rho}\,\nabla\cdot(\rho\,\textbf{a}\nabla)$. Intuitively, $X$ behaves like a diffusion process associated to the second order elliptic differential operator $\A$ in the interior of $D$, and is instantaneously reflected at the boundary in the inward {\it conormal} direction
$\displaystyle
\vec{\nu} =\textbf{a}\vec{n}$,
where $\vec{n}$ is the unit inward normal vector field on $\partial D$. See Chen \cite{zqC93} for the Skorokhod representation for $X$, which tells us some precise pathwise properties of $X$. We call $X$ an $(\textbf{a},\,\rho)$-\textbf{reflected diffusion} or an $(\A,\,\rho)$-\textbf{reflected diffusion}. A special but very important case is when $\textbf{a}$ is the identity matrix and $\rho=1$, in which $X$ is called a \textbf{reflected Brownian motion (RBM)}. Next, we make the following assumption about the killing potential throughout this paper.

\begin{assumption}\label{A:KillingPotential}(Killing potential)
Suppose $q(t,x)$ is a given non-negative bounded function on $[0,\infty)\times \bar{D}$ such that $q(t,\cdot)\in C(\bar{D})$ for all $t\geq 0$. Suppose also that $\delta_N$ is a sequence of positive numbers which converges to zero and denote
$q_N(t,x) = \delta_N^{-1}\,\1_{D^{\delta_N}}(x)q(t,x)$,  where
$D^{\delta} =\{x\in D:\; dist(x,\,\partial D)<\delta \}$.
\end{assumption}

Our particle system is parameterized by $N\in \mathbb{N}$, the initial number of particles. The function $q_N$ plays the role of a time-dependent killing potential. This killing potential is singular in the sense that $\delta_N^{-1}\,\1_{D^{\delta_N}}(x)$  converges weakly to the surface measure $\sigma$ which is singular with respect to Lebesque measure. More precisely, for $N\in \mathbb{N}$, we let $\{X_{i}\}_{i=1}^N$ be independent $(\textbf{a},\,\rho)$-reflected diffusions in $D$ and $\{R_i\}_{i=1}^N$ be independent exponential random variables with mean one. The normalized empirical measure of the particles \emph{alive} is defined as:
\begin{equation}\label{Def:EmpiricalRobin}
    \X^N_t(dz)  := \dfrac{1}{N}\sum_{\{i:\,t<\zeta^{(N)}_i \}}\1_{X_{i}(t)}(dz), \text{ where}
\end{equation}
\begin{equation}\label{E:KillingTime_N}
\zeta^{(N)}_i  = \inf \left\{t\geq 0: \frac{1}{2}\,\int_0^tq_N(s,X_i(s))ds  \geq R_i \right\}.
\end{equation}
Note that $\X^N_t$ is a random measure on $\bar{D}$. Moreover, $\X^N= (\X^N_t)_{t\geq 0}$ is a strong Markov process in $M_+(\bar{D})$, the space of finite non-negative Borel measures on $\bar{D}$ equipped with weak topology, and $\X^N$ has sample paths in the Skorokhod space $D([0,\,\infty),\,M_{+}(\bar{D}))$ almost surely.

\begin{remark}\label{Rk:sub-processes} \rm 
Let $\{Z^{(N)}_i\}_{i=1}^N$ be independent \textbf{sub-processes} (cf. \cite{CF12}) of reflected diffusions killed by the potential $q_N$. That is,
\[ Z^{(N)}_i(t) :=
    \begin{cases}
        X_i(t), \,& t<\zeta^{(N)}_i \\
        \partial, \,& t\geq \zeta^{(N)}_i,
    \end{cases}
\]
where $\partial$ is an isolated point of $\bar{D}$. Then $\X^N_t(dz)$ defined in \eqref{Def:EmpiricalRobin} is equal to $\frac{1}{N}\sum_{i=1}^{N}\1_{Z^{(N)}_{i}(t)}(dz)$ if we view $\1_{\partial}$ as the zero measure.
\end{remark}

We coin this model the name \textbf{Robin boundary model} due to the following hydrodynamic result. In what follows, $\toL$ denotes convergence in law and $\equalL$ denotes equal in law.

\begin{thm}\label{T:HydrodynamicLimit_RBM_Killing}\textbf{(Functional Law of Large Numbers)}
    Suppose Assumptions  \ref{A:GeometricSetting}  and  \ref{A:KillingPotential} hold. Suppose $\{\X^N_0\} \toL u_0(x)\rho(x)\,dx$ in $M_+(\bar{D})$, where $u_0\in C(\bar{D})$. Then
    $$\X^N_t(dx) \,\toL\, u(t,x)\rho(x)\,dx \quad \text{in  }D([0,\,\infty),\,M_{+}(\bar{D})),$$
    where $u\in C([0,\infty)\times\bar{D})$ is the probabilistic solution
    \footnote{
    By \cite{zqCwtF14b}, $u$ has the probabilistic representation $\E^{x}\left[u_0(X_t) \exp \left(-\int_0^tq(t-s,X_s)\,dL_s \right)\right]$ where $L_t$ is the boundary local time of $X$.
    }
    to the heat equation $\frac{\partial u}{\partial t}=\A u$ with Robin boundary condition $\frac{\partial u}{\partial \vec{\nu}}(t,x)=q(t,x)u(t,x)/\rho(x)$ on $(0,\infty)\times\partial D$   and initial condition $u(0,\cdot)=u_0$.
\end{thm}

The proof of Theorem \ref{T:HydrodynamicLimit_RBM_Killing} is an elementary law of large numbers argument involving the calculation of two moments. Since it is much easier than that of \cite{zqCwtF14b}, we omit it here and refer the reader to that paper.

\subsection{Main result}

Our object of study in this paper is the \textbf{fluctuation process} $\Y^N =(\Y^N_t)_{t\geq 0}$ defined by
\begin{equation}\label{e:1.4}
    \<\Y^N_t,\phi\> :=  N^{1/2}(\<\X^N_t,\phi\>-\E\<\X^N_t,\phi\>) \quad t\geq 0,\,\phi\in L^2(D),
\end{equation}
where $\<\X^N_t,\phi\> := \frac{1}{N}\sum_{\{i:\,t<\zeta^{(N)}_i \}}\phi(X_{i}(t))$ is the integral of an observable (or test function) $\phi$ with respect to  the measure $\X^N_t$. Even in this simple setting, answers to the following natural questions are non-trivial.
\begin{enumerate}
\item[(i)]   What is the state space for $\Y^N_t$? This space should posses a topology which allows us to make sense of convergence of $\Y^N$, if it does converge. Observe that although $\Y^N$ acts on $L^2(D)$ linearly, it is not a bounded operator in general.
\item[(ii)]   Does $\Y^N$ converge? If so, what can we say about its limit?
\end{enumerate}
The answer for question (i) is given by Lemma \ref{L:StateSpace_RBMkilling_Y}. It says that the process $(\Y^N_t)_{t\geq 0}$ has sample paths in $D([0,\infty),\H_{-\alpha})$ for $\alpha>0$ large enough, where $\H_{-\alpha}$ is a Hilbert space of distributions that strictly contains $L^2(D,\,\rho(x)dx)$. See subsection \ref{subsubsection:HilbertSpaces} for the precise construction of $\H_{-\alpha}$, which can be identified with the dual of the Sobolev space
$W^{\alpha/2, 2}(D)$  of fractional order.

The answer for question (ii) is given by Theorem \ref{T:Convergence_RBM_killing_Y},
the main result of this paper. Theorem \ref{T:Convergence_RBM_killing_Y} contains 2 parts: the convergence result and the properties of the limit. The limit is shown to be decomposable into an independent sum of a ``transportation part" and a ``white noise part" (see \eqref{E:GeneralizedOUformula_RBM_Killing}  below). The `transportation part' is governed by the evolution operators $\{Q_{s,t}\}_{s\leq t}$ generated on $C(\bar{D})$ by the backward PDE $\frac{\partial v}{\partial s}=- \A v$ on $(0,t)\times D$ with Robin boundary condition $\frac{\partial v}{\partial \vec{n}} = qv/\rho$ on $(0,t)\times \partial D$. More precisely, for $0 \leq s\leq t$ and $\phi\in L^2(D)$, we define
\begin{eqnarray}\label{E:PbtyRepresentation_Q}
Q_{s,t}\phi(x)
&:=&  \E\left[\phi(X_{t })\, \exp \left(-\int_s^t q(r,X_r)\,dL_r \right) \Big|\,X_s=x\right] \nonumber \\
& = &\E\left[\phi(X_{t-s})\, \exp \left(-\int_0^{t-s} q(s+r,X_r)\,dL_r \right) \Big|\,X_0=x\right] .
 \end{eqnarray}
Define
\begin{equation}\label{Def:U_ts_RBMkill}
\mathbf{U}_{(t,s)}\mu (\phi):= \mu( Q_{s,t}\phi)
\end{equation}
for $\alpha>0$, $\mu\in \H_{-\alpha}$ and $\phi\in L^2(D)$ whenever it is well defined (i.e. $Q_{s,t}\phi\in \H_{\alpha}$); see
Theorem \ref{T:Convergence_RBM_killing_Y} and Remark \ref{R:2.1}.
For simplicity,   denote by $\<\phi,\,\psi\>_{\rho}:= \int_D \phi(x)\psi(x)\,\rho(x)dx$ the inner product of $L^2(D,\,\rho(x)dx)$. We can now formulate our main result.

\begin{thm}\label{T:Convergence_RBM_killing_Y}\textbf{(Functional Central Limit Theorem)}
    Suppose that Assumptions \ref{A:GeometricSetting}  and  \ref{A:KillingPotential}  hold and that the initial positions of particles are i.i.d with distribution $u_0(x)\rho(x)\,dx$, where $u_0\in C(\bar{D})$. Then for any $\alpha>d+2$ and $T>0$, $\Y^N$ converges to $\Y$ in distribution as $N\to \infty$ in the Skorokhod space $D([0,\,T],\H_{-\alpha})$, where $\Y$ is the generalized Ornstein-Uhlenbeck process taking values in $D([0,\,T],\H_{-\alpha})$ given by
    \begin{equation}\label{E:GeneralizedOUformula_RBM_Killing}
    \Y_t \equalL \mathbf{U}_{(t,0)}\Y_0 + \int_0^t \mathbf{U}_{(t,s)}\,dM_s .
    \end{equation}
    In the above, $M$ is a (unique in distribution) continuous, $\tilde{\mathcal{F}}_{t}$-adapted, square integrable, $\H_{-\alpha}$-valued Gaussian martingale with independent increments and covariance functional characterized by
    \begin{equation}\label{E:covariance_M_Killling}
    \tilde{\E} \left[\<M_t,\phi\>^2 \right]
    = \int_0^t \bigg( \<\textbf{a}\nabla \phi \cdot \nabla \phi,\,u(s)\>_{\rho} + \int_{\partial D}\phi^2(z)u(s,z)q(s,z)\rho(z)\,d\sigma(z) \bigg)\;ds,
    \quad \phi \in \H_\alpha,
    \end{equation}
    defined on a complete probability space with right continuous filtration $(\tilde{\Omega},\,\tilde{\mathcal{F}},\,\tilde{\mathcal{F}}_{t},\,\tilde{\P})$, where the function $u(s,x)$ is given by Theorem \ref{T:HydrodynamicLimit_RBM_Killing}. $\Y_0$ is the centered Gaussian random variable with covariance
    $$
    \tilde{\E}\,[\Y_0(\phi)\Y_0(\psi)]= \<\phi\,\psi,\,u_0\>_{\rho}- \<\phi,u_0\>_{\rho}\,\<\psi,u_0\>_{\rho}
     \quad \hbox{for } \phi , \psi \in \H_\alpha,
    $$
    defined on the same probability space as $M$ and is independent of $M$.
    Moreover, $\Y$ is a continuous Gaussian Markov process which is unique  in distribution, and $\Y$ has a version in $C^{\gamma}([0,\infty),\,\H_{-\alpha})$ (i.e. H\"older continuous with exponent $\gamma$) for any $\gamma\in(0,1/2)$.
\end{thm}

\begin{remark}\label{R:1.6}\rm \begin{description}
\item{(i)} In \eqref{E:GeneralizedOUformula_RBM_Killing}, $\int_0^t \mathbf{U}_{(t,s)}\,dM_s$ is the stochastic integral with respect to  the Hilbert space valued martingale $M$ (cf. \cite{MP80}). In the Appendix, we prove that it is well-defined. For the convenience of the reader, we also stated the precise definition of Hilbert space valued continuous Gaussian processes with independent increment. The existence and uniqueness of $M$ is given in Theorem \ref{L:ConvergenceOfM^N_RBM_killing}. Furthermore, for $\alpha>d+2$, both $\mathbf{U}_{(t,0)}\Y_0$ and $\int_0^t \mathbf{U}_{(t,s)}\,dM_s$ live in $\H_{-\alpha}$ (i.e. they extend to be continuous functionals on $\H_{\alpha}$).

\item{(ii)} Roughly speaking, $\Y$ solves the following stochastic evolution equation (called the \textbf{Langevin equation}) in the weak sense:
\begin{equation}\label{E:SEE_RBM_killing_limit}
    dY_t = \mathbf{A}^{(-\alpha)}_tY_t\,dt + dM_t, \quad Y_0=\Y_0,
\end{equation}
where $\mathbf{A}^{(-\alpha)}_t$ is the generator of $\{\mathbf{U}_{(t,s)}\}_{t\geq s}$ in the Hilbert space $\H_{-\alpha}$.

\item{(iii)}
Define a bilinear forms $\D^{(q)}_s$ on $L^2(D,\,\rho(x)dx)\cap L^2(\partial D,d\sigma)$ by
    \begin{equation}\label{Def:D_r(phipsi)}
    \D^{(q)}_s(\phi,\psi):= \<\textbf{a}\nabla \phi \cdot \nabla \psi,\,u(s)\>_{\rho} + \int_{\partial D}\phi\,\psi \,u(s)\,q(s)\,\rho\,d\sigma
    \end{equation}
and $\D^{(q)}_s(\phi):=\D^{(q)}_s(\phi,\phi)$ for $s\geq 0$. Now \eqref{E:covariance_M_Killling}  reads as $\tilde{\E}[\<M_t,\phi\>^2] = \int_0^t \D^{(q)}_s(\phi)\,ds$. As an immediate application of \eqref{E:GeneralizedOUformula_RBM_Killing}, for all fixed $\phi\in \H_{\alpha}$ with $\alpha>d+2$, we have
    \begin{equation}\label{E:GeneralizedOUformula_RBM_Killing_eg}
    \Y_t(\phi) \equalL \Y_0(Q_{0,t}\phi) + \int_0^t \sqrt{\D^{(q)}_{s}(Q_{s,t}\phi)}\,dB^{(\phi)}_s \quad \text{ in } D([0,\,T],\R),
    \end{equation}
where $B^{(\phi)}$ is a standard Brownian motion independent of $\Y_0$. Therefore, we can simulate the evolution (in time $t$) of the fluctuations of the particle density with respect to  an observable $\phi$
by running a Brownian motion.

\item{(iv)} When $D$ is a cube (such as when $d=1$), Theorem \ref{T:Convergence_RBM_killing_Y} holds with $\alpha>d/2+2$ in place of $\alpha>d+2$, since we have a stronger unform upper bound for eigenfunctions, namely $\sup_{\ell}\|\phi_{\ell}\|<C(d,D)$. 
    \qed
\end{description}
\end{remark}

\begin{remark}\label{R:1.7}\rm \begin{description}
\item{(i)} When $q=0$,  Theorem \ref{T:Convergence_RBM_killing_Y} in particular
gives the fluctuation result for independent reflecting Brownian motions in bounded Lipschitz domains.

\item{(ii)} (Killing by local time) \rm
 Clearly, the measure $q_N(t,x)\,dx$  converges weakly to $q(t,x)d \sigma(x)$ as $N\to\infty$, where $\sigma$ denotes the surface measure on $\partial D$. The positive additive continuous functional (see the Appendix of \cite{CF12}) of $X_i$ having Revuz measure $q(t,x)d \sigma(x)$ is $2\int_0^tq(s,X_i(s))dL^{(i)}_s$, where $L^{(i)}_t$ is the boundary local time of $X_i$. Hence it is natural to ask:  what if the processes $\{X_i\}_{i\geq 1}$ are killed by $2\int_0^tq(s,X_i(s))dL^{(i)}_s$ (which no longer depends on $N$) rather than by a potential function $q_N$ on the strip $D^{\delta_N}$? It turns out that, with little extra effort, one can show that Theorem \ref{T:HydrodynamicLimit_RBM_Killing} and Theorem \ref{T:Convergence_RBM_killing_Y} both remain valid if we replace the definition of $\zeta^{(N)}_i$ in \eqref{E:KillingTime_N} by
\begin{equation}\label{E:KillingTime}
\zeta^{(N)}_i= \zeta_i:= \inf \left\{t\geq 0:\; 2\int_0^tq(s,X_i(s))dL^{(i)}_s  \geq R_i \right\}.
\end{equation}
See subsection \ref{S:5.4}  for details.
\qed
\end{description}
\end{remark}

One of the earliest rigorous results about fluctuation limit was proven by It\^{o}   \cite{kI80, kI83}, who considered a system of independent and identically distributed (i.i.d.) Brownian motions in $\R^d$ and showed that the limit is a $\mathcal{S}'$-valued Gaussian process solving a Langevin equation, where $\mathcal{S}'$ is the Schwartz space of tempered distributions. Fluctuation limits for stochastic particle systems in domains are very limited. Sznitman \cite{asS84} studied the fluctuations of a conservative system of diffusions with normal reflected boundary conditions on smooth domains. Fluctuations of the reaction-diffusion systems on the cube $[0,1]^d$ with linear or quadratic reaction terms were studied in \cite{BDP92,pD88b,pK86, pK88}. These fluctuation results are valid only for dimension $d\leq 3$.

\subsection{Outline of proof}\label{S:1.2}

We prove Theorem \ref{T:Convergence_RBM_killing_Y} through the following six steps.

\begin{enumerate}
\item[] Step 1: \; $\Y^N$ satisfies the following stochastic integral equation
                \begin{eqnarray*}
                    \Y^{N}_t= \mathbf{U}^N_{(t,0)}\Y^{N}_0+ \int_0^t \mathbf{U}^N_{(t,s)}\,dM^N_s\;\text{a.s.},
                \end{eqnarray*}
                \qquad \qquad where $\mathbf{U}^N_{(t,s)}$ is an evolution system approximating $\mathbf{U}_{(t,s)}$; see Theorem \ref{T:3.3}.
\item[]  Step 2: \; $M^N \toL M$ in $D([0,T],\H_{-\alpha})$; see Theorem \ref{L:ConvergenceOfM^N_RBM_killing}.
\item[]  Step 3:\; $\Y^N$ is tight in $D([0,T],\H_{-\alpha})$; see Theorem  \ref{T:3.7}.
\item[]  Step 4: \; $\mathbf{U}^N_{(t,0)}\Y^N_0 \toL \mathbf{U}_{(t,0)}\Y_0$ in $D([0,T],\H_{-\alpha})$; see Theorem \ref{T:3.8}.
\item[] Step 5: \; $\int_0^t \mathbf{U}^N_{(t,s)}\,dM^N_s \toL  \int_0^t \mathbf{U}_{(t,s)}\,dM_s$ in $D([0,T],\H_{-\alpha})$; see Theorem \ref{T:ConvergenceStocInt}.
\item[] Step 6: \; All the stated properties for the fluctuation limit hold; see Theorem \ref{T:Convergence_RBM_Y}.
\end{enumerate}

The main difficulty is in establishing the convergence in Step 5.
Note that $t\mapsto \int_0^t \mathbf{U}_{(t,s)}\,dM_s$ is not a martingale. The standard method based on Kotelenez's submartingale inequality \cite{pK82} does not seem to work. This is because in our case $\mathbf{U}_{(t,s)}$ is not exponentially bounded; that is, there is no $\beta >0$ so that the operator norm $\|\mathbf{U}_{(t,s)}\| \leq e^{\beta (t-s)}$ for $t\geq s$ (see \cite{pK82}). In fact,  we suspect it is not even a bounded operator on $\H_{-\alpha}$ due to the singular interaction near the boundary. To overcome this difficulty, we need first to make sense of the expression $\int_0^t \mathbf{U}_{(t,s)}\,dM_s$, which is done in Section 4, the Appendix of this paper. Our approach is then based on  suitably extending the functional analytic framework of \cite{pK86} and   a direct analysis that uses heat kernel estimates and Dirichlet Form method.

\section{Functional analytic framework}

Our method to study the fluctuation is functional analytic, with the mathematical framework being the calculus of evolution equations on Hilbert spaces (see, for example, \cite{pDjZ08, GT95, mH09}). As remarked in \cite{pK86}, this approach yields a useful representation of the limiting process (the generalized Ornstein-Uhlenbeck process) as the mild solution of a stochastic partial differential equation (SPDE), which yields uniqueness and Gaussian property for free. It also tells us the smallest Hilbert space in which the generalized Ornstein-Uhlenbeck process lives.

\textbf{Conventions and notations: }

In this paper, we use := as a way of definition. For $a, b\in \R$, $a\vee b:= \max \{a, b\}$ and $a\wedge b := \min \{a, b\}$. We use abbreviation r.c.l.l.  for right continuous having left limits, and $\|\cdot\|$ to denote the supremum norm in $\bar{D}$. Even though the constants appearing in the article may depend on $\textbf{a}$ or $\rho$ given in Assumption \ref{A:GeometricSetting}, we will not mention this dependence explicitly. For example, we use $C(d,D)$ to denote a constant which depends only on $d$ and $D$ (and possibly on $\textbf{a}$ or $\rho$). The exact value of the constant may vary from line to line.

\subsection{Neumann heat kernel}

It is well known (cf. \cite{BH91, GSC11} and the references therein) that, on a bounded Lipschitz domain $D$, an $(\A,\,\rho)$-reflected diffusion $X$ has a jointly locally H\"older continuous transition density $p(t,x,y)$ with respect to the symmetrizing measure $\rho(x)dx$ on $(0,\infty)\times \bar{D}\times \bar{D}$. Moreover, the following Aronson type Gaussian estimates hold:
   \begin{equation}\label{e:2.1}
        \dfrac{1}{c_1 t^{d/2}}\,\exp\left(\frac{-c_2 |y-x|^2}{t}\right)
\leq p(t,x,y) \leq \dfrac{c_1}{t^{d/2}}\,\exp\left(\frac{-|y-x|^2}{c_2\,t}\right)
    \end{equation}
for every $(t,x,y)\in(0,T]\times \bar{D}\times \bar{D}$, where $c_1=c_1(d,D,T)$ and $c_2=c_2(d,D,T)$ are positive constants.

Using \eqref{e:2.1} and the Lipschitz assumption for $\partial D$, we can check that
\begin{equation}\label{E:BoundEpsilonStrip}
\sup_{\eps\in(0,\eps_0)}\sup_{x\in \bar{D}}\frac{1}{\eps}\int_{D^{\eps}}p(t,x,y)\,dy \leq \frac{C}{t^{1/2}}\quad \text{for }t\in(0,T],
\end{equation}
where $\eps_0=\eps_0(D)>0$ and $C=C(d,D,T)>0$. In particular, we can let $\eps\to 0$ in \eqref{E:BoundEpsilonStrip} to obtain, via \eqref{L:MinkowskiContent_D},
\begin{equation}\label{E:BoundPartialD}
\sup_{x\in \bar{D}}\int_{\partial D}p(t,x,y)\,d\sigma(y) \leq \frac{C}{t^{1/2}}\quad \text{for }t\in(0,T].
\end{equation}

\subsection{Hilbert space  $\H_{\gamma}$}\label{subsubsection:HilbertSpaces}

Recall that $\A = \frac{1}{2\,\rho}\,\nabla\cdot(\rho\,\textbf{a}\nabla)$ denotes the $L^2(D,\,\rho(x)dx)$-generator for an $(\textbf{a},\,\rho)$-reflected diffusion. Clearly, $\A$ is a self-adjoint, non-positive operator on $L^2(D,\,\rho(x)dx)$. Together with the fact that $D$ is bounded, we see that $\A$ has a discrete spectrum in $\H_0$. Let $\phi_{k}$ be a complete orthonormal system (CONS) of eigenvectors of $\A$ in $\H_0$ with eigenvalues $-\lambda_{k}$, where $0<\lambda_{1}\leq \lambda_2 \leq \lambda_3 \leq \cdots$. Note that the linear span of $\{\phi_{k}\}$ is dense in $L^2(D; \rho dx)$. We define, for $\alpha\in (-\infty,\infty)$,
\begin{equation}
\H_{\alpha}:= \text{ the closure of the linear span of }\{\phi_{k}\}  \text{ with respect to  the inner product }\<\,,\,\>_{\alpha},
\end{equation}
where $\<\phi,\psi\>_{\alpha}:= \<(I-\A)^{\alpha}\phi,\,\psi\>_{\rho}$. Here $I$ is the identity operator on $\H_0=L^2( D; \rho dx)$ and $(I-\A)^{\alpha}$ is the $\alpha$-th power (defined through spectral representation) of the positive self-adjoint operator $I-\A$. In particular, $\<\,,\,\>_{0}=\<\,,\,\>_{\rho}$ by definition.

Note that $(\H_{\alpha},\,\<\,,\,\>_{\alpha})$ is a real separable Hilbert space and that $\H_{\beta}\subset \H_{\alpha}$ when $\beta>\alpha$. Moreover, $\H_\alpha$ and  $\H_{-\alpha}$ are dual to each other.
Equip $\Phi := \cap_{\alpha \geq 0} \H_{\alpha}$ with the locally convex topology defined by the set of norms $\{|\varphi|_{\alpha}:= \<\varphi,\varphi\>_{\alpha}^{1/2}:\;\varphi\in\Phi,\,\alpha\in[0,\infty)\,\}$. Let $\Phi'$ be the strong dual of $\Phi$. Identifying $\H_0$ with its dual $\H_0'$, we obtain the chain of dense continuous inclusions
\begin{equation}\label{e:2.4}
\Phi\subset \H_{\alpha} \subset \H_0=\H_0' \subset \H_{-\alpha} \subset \Phi', \quad \alpha\in[0,\infty).
\end{equation}
Moreover, for $\beta\in \R$, we have
\begin{equation}\label{e:2.5}
h^{(\beta)}_{k}:= (1+\lambda_{k})^{-\beta/2}\,\phi_{k} \quad \text{is a CONS for }\H_{\beta}.
\end{equation}
Hence, $\<\phi,\psi\>_{\beta}=\sum_{k\geq 1}\,\<\phi\,,\,h^{(-\beta)}_{k}\,\>\,\<\psi\,,\,h^{(-\beta)}_{k}\,\>$ for $\phi,\,\psi\in \H_{\beta}$ and
\begin{equation}
\H_{\beta}=\Big\{\mu\in\Phi':\;\sum_{k\geq 1}\<\mu\,,\,h^{(-\beta)}_{k}\,\>^2\,<\infty \Big\},
\end{equation}
where $\<\,,\,\>$ denotes the dual paring extending $\<\,,\,\>_{\rho}$.

\begin{remark}\label{R:2.1} \rm
When $\alpha >0$, $\H_\alpha$ can be identified with the fractional Sobolev
space $W^{\alpha/2, 2}(D)$ on $D$. This is because for
$\alpha \geq 0$, $\H_\alpha = (I -\A)^{-\alpha/2} L^2(D,\,\rho dx)$. Since $D$ is a bounded Lipschitz domain,
it is known that $\H_\alpha = W^{\alpha/2, 2}(D)$ when $\alpha =1$
(see \cite{zqC93}) and hence for every integer $\alpha\geq 1$.
It follows by interpolation that $\H_\alpha = W^{\alpha/2, 2}(D)$
for every $\alpha >0$. When $\alpha<0$, $\H_\alpha $ can be identified as the
dual space of $\H_{-\alpha}$. \qed
\end{remark}

\subsection{Weyl's law and eigenfunction estimates}

For a general bounded Lipschitz domain $D\subset \R^d$, the Weyl's asymptotic law for the Neumann eigenvalues   holds (see \cite{NS05}). That is, the number of eigenvalues (counting their multiplicities) less than or equal to $x$, denoted by $\sharp\,\{k:\,\lambda_{k}\leq x\}$, satisfies
\begin{equation}\label{E:WeylLaw}
\lim_{x\to\infty}\frac{\sharp\,\{k:\,\lambda_{k}\leq x\}}{x^{d/2}}=C \quad \text{ for some constant }C=C(d,D)>0.
\end{equation}

\begin{lem}
There exists $C=C(d,D)>0$ such that for all integers $k\geq 1$ we have
\begin{equation}\label{E:EigenfcnUpperBound}
\|\phi_{k}\|  \leq C\,\lambda_{k}^{d/4} \quad \hbox{ and }
\quad \int_{\partial D} \phi_{k}^2\,d\sigma \leq C\,(\lambda_{k}+1).
\end{equation}
\end{lem}

\begin{pf}
By Cauchy-Schwartz inequality, Chapman Kolmogorov equation and then the Gaussian upper bound, we have
\begin{eqnarray*}
|\phi_{k}(x)| &=& e^{\lambda_{k}t}|P_{t}\phi_{k}(x)| \leq e^{\lambda_{k}t}\|\phi_{k}\|_{L^2(\rho)}\,\sqrt{p(2t,x,x)}\\
&\leq& e^{\lambda_{k}t}\,\frac{C(d,D)}{t^{d/4}} \quad\text{for }t\leq 1/\lambda_1.
\end{eqnarray*}
Taking $t=1/\lambda_{k}$ yields the first inequality in \eqref{E:EigenfcnUpperBound}.

Recall that the Dirichlet form $(\D,\,Dom(\D))$ (in $L^2(D,\,\rho(x)dx)$) for the $(\A,\,\rho)$-reflected diffusion $X$ is regular (since $D$ has Lipschitz boundary (cf. \cite{BH91})) and that the surface measure $\sigma$ is smooth. Hence by Theorem 2.1 of \cite{SV96}, we have the following generalized trace theorem:
\begin{equation}
\int_{\partial D} f(x)^2\,\sigma(dx) \leq \| G_\beta \sigma \|\left(\,\D(f, f)+ \beta \int_{D} f^2(x)\,dx\,\right)
\end{equation}
for any $f\in Dom(\D)$ and $\beta>0$, where $G_\beta \sigma (x):= \int_0^\infty\int_{\partial D} e^{-\beta t} p(t, x, y)\,\sigma(dy)dt$. Note that $\|G_\beta \sigma \| <\infty$ by \eqref{E:BoundPartialD} and the fact that $p(t,x,y)$ converges to $1/\int_D\rho (x) dx$  as $t\to\infty$  uniformly for $(x,y)\in\bar{D}\times \bar{D}$ exponentially fast (by eigenfunction expansion). Hence, taking $\beta=1$, we obtain the second inequality in \eqref{E:EigenfcnUpperBound}.
\end{pf}

\section{Preliminaries}

\subsection{Minkowski content for $\partial D$}

By the same proof of \cite[Lemma 7.1]{zqCwtF14b}, we obtain the following result.
\begin{lem}\label{L:MinkowskiContent_D}
    Let $D\subset \R^d$ be a bounded Lipschitz domain and $k\in\mathbb{N}$. If $\F\subset  C(\bar{D}^k)$ is an equi-continuous and uniformly bounded family of functions, then
    $$ \lim_{\eps \to 0}
    \frac{1}{\eps^k}\int_{(D^{\eps})^k}f(z_1,\cdots,z_k)\,dz_1\cdots dz_k = \int_{(\partial D)^k}f(z_1,\cdots,z_k)\,\sigma(dz_1)\cdots \sigma(dz_1)$$
    uniformly for $f\in\F$, where $D^{\eps}:=\{x\in D:\; dist(x,\,\partial D)<\eps \}$ is the $\eps$-neighborhood of $\partial D$ in $D$ and $\sigma$ is the surface measure on $\partial D$.
\end{lem}

By a simple modification of the same proof, we can strengthen the above lemma as follows.

\begin{lem}\label{L:MinkowskiContent_D_I}
    Let $D\subset \R^d$ be a bounded Lipschitz domain, $I$ be a $\mathcal{H}^{d-1}$-rectifiable closed subset of $\partial D$ and $k\in\mathbb{N}$. If $\F \subset \mathcal{B}(D^k)$ is an equi-continuous and uniformly bounded family of functions on an open neighborhood of $I^k$, then
    $$
    \lim_{\eps \to 0} \frac{1}{\eps^k}\int_{(I^{\eps})^k}f(z_1,\cdots,z_k)\,dz_1\cdots dz_k =  \int_{I^k}f(z_1,\cdots,z_k)\,\sigma(dz_1)\cdots \sigma(dz_1)$$
    uniformly for $f\in\F$, where $I^{\eps}:=\{x\in D:\; dist(x,\,I)<\eps \}$  is the $\eps$-neighborhood of $I$ in $D$.
\end{lem}

The following is about a convergence result uniform in the shrinking rate of $\delta=\delta_N$. It is used to guarantee that $\delta_N$ can be any sequence (which converges to zero) in the proof of Lemma \ref{L:ConvergenceDistOfIntegral}.

\begin{lem}\label{L:HNGNtoZero}
Suppose $\{\X^N_0\} \toL u_0(x)\rho(x)\,dx$ in $M_+(\bar{D})$ as $N\to \infty$, where $u_0\in C(\bar{D})$.
Let $\{\varphi_N(r)\,:\,r\geq 0,\,N\in\mathbb{N}\}$ be a family of non-negative continuous functions on $\bar{D}$
such that $\sup_{N}\sup_{r\geq 0}\|\varphi_N(r)\| <\infty$. For any $\delta_N\to 0$, $T>0$ and $p\geq 1$, we have
\begin{equation}\label{E:HNGNtoZero}
\lim_{N\to\infty} \E\left[\,\bigg(\sup_{t\in[0,T]}\Big| \int_0^t\<\varphi_N(r) \delta_N^{-1}\1_{D^{\delta_N}},\,\X^N_r\>
- \< \varphi_N(r)  \delta_N^{-1} \1_{D^{\delta_N}},\,u(r)\>_{\rho}\,dr\Big|\bigg)^p\,\right] =0,
\end{equation}
where $\1_{D^{\delta_N}}$ is the indicator function on $D^{\delta_N}$.
\end{lem}

\begin{pf}
Let
$H_N(t):=\int_0^t\< \varphi_N(r) \delta_N^{-1} \1_{D^{\delta_N}},\,\X^N_r\>\,dr $ and
$G_N(t):=\int_0^t \< \varphi_N(r) \delta_N^{-1}\1_{D^{\delta_N}},\,u(r)\>_{\rho}\,dr$.
It can be shown, by a standard argument and using Lemma \ref{L:MinkowskiContent_D}, that for any $T>0$,
\begin{equation*}
H_N(t)-G_N(t) \,\toL \,0 \quad \text{in } \,C([0,T],\R).
\end{equation*}
In particular, by the metric of $C([0,T],\R)$ and the deterministic nature of the limit, we have
\begin{equation*}
\sup_{t\in[0,T]}|H_N(t)-G_N(t)| \to 0 \quad \text{both in law and in probability}.
\end{equation*}
On other hand, since $H_N(t)$ and $G_N(t)$ are increasing, we have
\begin{equation*}
\limsup_{N\to\infty}\,\E\left[\,\Big(\sup_{t\in[0,T]}|H_N(t)-G_N(t)|\Big)^p\,\right] \leq \limsup_{N\to\infty}\,2^p\,\Big(\E\left[\,H^p_N(T)\,\right] + G^p_N(T)\Big).
\end{equation*}
Furthermore, we can check that $\limsup_{N\to\infty}\E[H^p_N(T)+ G^p_N(T)] <\infty$.
Denote by ${\cal P}(\overline D) $ the collection of sub-probability measures on $\overline D$.
Comparing with the process without killing (i.e. replacing the sub-processes $Z^{(N)}_i$ by the reflected diffusions $X_i$ in the definition of $\X^N$),
we have by \eqref{E:BoundEpsilonStrip}
\begin{eqnarray*}
\sup_{\mu \in {\cal P}(\overline D)} \E_\mu \left[ H_N(t)\right]
\leq \| \varphi_N\|
\sup_{x\in \overline D} \E_x  \int_0^t  \1_{D^{\delta_N}} (X_1 (r)) dr
= \| \varphi_N\| \sup_{x\in \overline D}   \int_0^t  \int_{D^{\delta_N}} p(r, x, y) dy dr \leq C_1\, t^{1/2},
\end{eqnarray*}
where $C_1$ is a positive constant independent of $N$ and $t$.
Let $f(r):=\< \delta_N^{-1}\1_{D^{\delta_N}},\,\X^N_{r}\>$.
Then for any positive integer $k$, by Fubinni's theorem and the Markov property, we have for any initial distribution $\mu$ of $\X^N_{0}$,
\begin{equation*}
\E_\mu [H^k_N(T)]= k!\,\E \int_{0\leq r_1 \leq r_2 \leq \cdots \leq r_k\leq T}
 f(r_1) f(r_2) \dots f (r_k) dr_1 dr_2 \cdots dr_k \leq k! \, ( C_1 T^{1/2})^k.
\end{equation*}
It in particular implies that, under the assumption  $\{\X^N_0\} \toL u_0(x)\rho(x)\,dx$ in $M_+(\bar{D})$,
$$
\limsup_{N\to \infty} \E[H^k_N(T)] \leq  \| u_0\| \| \rho \| \,
k! \, ( C_1 T^{1/2})^k.
$$
A similar argument yields  $\limsup_{N\to\infty}\E[G^k_N(T)] <\infty$ for any positive integer $k$. Hence, by interpolation, we have $\limsup_{N\to\infty}\E[H^p_N(T)+ G^p_N(T)] <\infty$ for all $p\geq 1$.

The uniform integrability implied by $\limsup_{N\to\infty}\,\E\left[\,\Big(\sup_{t\in[0,T]}|H_N(t)-G_N(t)|\Big)^p\,\right] <\infty$, together with the convergence $\sup_{t\in[0,T]}|H_N(t)-G_N(t)|\to 0$ in probability, guarantee (see, e.g. Theorem 5.2 in \cite[Chapter 4]{Durrett10}) that the lemma is true.
\end{pf}

\subsection{Estimates for evolution semigroups $Q^N_{(s,t)}$ and $Q_{(s,t)}$}

Recall the definition of $Q_{(s,t)}$ and $\mathbf{U}_{(t,s)}$ in \eqref{E:PbtyRepresentation_Q} and \eqref{Def:U_ts_RBMkill}, respectively.
For any fixed $t>0$ and $\phi\in C(\bar{D})$,
$v(s,x):=Q_{(s,t)}\phi(x)$ is the unique element in $C([0,t]\times \bar{D})$ satisfying the integral equation
\begin{equation}\label{E:IntegralEqt_Q}
v(s,x)=P_{t-s}\phi(x)-\dfrac{1}{2}\int_0^{t-s}\int_{\partial D}p(\theta,x,y)\,q(s+\theta,y)\,v(s+\theta,y)\,\rho(y)\,d\sigma(y)\,d\theta ;
\end{equation}
see \cite[Proposition 4.1]{zqCwtF14b}.
We call $v$ the \textbf{probabilistic solution} of the backward equation

    \begin{equation}\label{E:BackwardHeatEqt_Q}
        \left\{\begin{aligned}
        \dfrac{\partial v}{\partial s} &= -\A v  & &\qquad\text{on } (0,t)\times D  \\
        \dfrac{\partial v}{\partial \vec{\nu}} &= \frac{q\,v}{\rho} & &\qquad\text{on } (0,t)\times \partial D  \\
        v(t)&= \phi & &\qquad\text{on } D.
        \end{aligned}\right.
    \end{equation}

Analogous to the definition of $Q_{(s,t)}$ and $\mathbf{U}_{(t,s)}$, we define
\begin{eqnarray}\label{E:PbtyRepresentation_QN}
Q^{N}_{s,t}\phi(x) &:=& \E\left[\phi(X_{t })\,\exp \left(-\int_s^{t }q_N( r,X_r)\,dr\right) \Big|\,X_s=x\right] \nonumber \\
&=& \E\left[\phi(X_{t-s})\,\exp \left(-\int_0^{t-s}q_N(s+r,X_r)\,dr\right) \Big|\,X_0=x\right]
\end{eqnarray}
and
\begin{equation}\label{Def:UN_ts_RBMkill}
\mathbf{U}^N_{(t,s)}\mu (\phi):= \mu( Q^N_{s,t}\phi)
\end{equation}
for $\alpha>0$, $\mu\in \H_{-\alpha}$ and $\phi\in L^2(D)$ whenever it is well defined (i.e. $Q^N_{s,t}\phi\in \H_{\alpha}$).
Then $v_N(s,x):=Q^{N}_{(s,t)}\phi(x)$ is the unique element in $C([0,t]\times \bar{D})$ satisfying the integral equation
\begin{equation}\label{E:IntegralEqt_QN}
v_N(s,x)=P_{t-s}\phi(x)-\dfrac{1}{2}\int_0^{t-s} P_{\theta}\left(q_N(s+\theta)v_N(s+\theta)\right)(x)\,d\theta,
\quad 0\leq s\leq t,
\end{equation}
provided that $\phi\in C(\bar{D})$. Here $\{P_t; t\geq 0\}$ is the transition semigroup of $X$ in $L^2(D,\,\rho(x)dx)$ (i.e. $P_tf(x)=\E^x[f(X_t)]= \int_D f(y)p(t,x,y)\rho(y)\,dy$).
As before, $v_N$ is called the \textbf{probabilistic solution} of the backward equation
\begin{align}\label{E:BackwardHeatEqt_QN}
    \begin{cases}
        \dfrac{\partial v_N}{\partial s} = - \frac{1}{2}\Delta v_N + q_N\,v_N  &\text{ on  } (0,t)\times D  \smallskip \\
        \dfrac{\partial v_N}{\partial \vec{\nu}} = 0 &\text{ on  } (0,t)\times \partial D  \smallskip \\
        v_N(t)=\phi &\text{ on  } D
    \end{cases}
\end{align}

\begin{remark} \rm
It can be shown (cf. \cite{CF12}), using the Markov property of the reflected diffusion $X$, that each $Z=Z^i$ (described in Remark \ref{Rk:sub-processes}) is a time-inhomogeneous Markov process on $\bar{D}\cup \{\Delta^{(i)}\}$ with $(Q^{N}_{s,t})_{s\leq t}$ being its transition operator: $Q^N_{s,t}f(x)=\E[f(Z_t)|Z_s=x]$, with the convention that $f(\Delta)=0$. Besides, (\ref{E:BackwardHeatEqt_QN}) is the Kolmogorov's backward equation for $Z$ and (\ref{E:PbtyRepresentation_QN}) is the probabilistic representation of the solution to (\ref{E:BackwardHeatEqt_QN}).  \qed
\end{remark}

The following uniform convergence and uniform bound are useful in many places of this paper.

\begin{lem}
For all $\phi\in C(\bar{D})$ and $0\leq s\leq t$, we have
\begin{equation}\label{E:UniformConvergenceQNQ}
\lim_{N\to\infty}Q^N_{s,t}\phi = Q_{s,t}\phi \quad\text{uniformly on }\bar{D}  \quad\text{and}
\end{equation}
\begin{equation}\label{E:ContractionQNQ}
\sup_N |Q^N_{s,t}\phi(x)|\vee |Q_{s,t}\phi(x)|\leq P_{t-s}|\phi|(x) \leq \|\phi\| \quad\text{for }x\in \bar{D}.
\end{equation}
\end{lem}

\begin{pf} Estimates \eqref{E:ContractionQNQ} follows immediately from (\ref{E:PbtyRepresentation_QN}), (\ref{E:PbtyRepresentation_Q}) and the non-negativity of $q$.  For (\ref{E:UniformConvergenceQNQ}), note that
\begin{eqnarray*}
&& \Big|Q^N_{s,t}\phi(x) - Q_{s,t}\phi(x)\Big| = \Big|\E_x\left[\phi(X_{t-s})\,\left(e^{-\int_0^{t-s}q_N(s+r,X_r)\,dr}- e^{-\int_0^{t-s}q(s+r,X_r)\,dL_r}\right)\right]\Big|\\
&\leq& \|\phi\|\,\E_x\left[\Big|\int_0^{t-s}q_N(s+r,X_r)\,dr- \int_0^{t-s}q(s+r,X_r)\,dL_r\Big|^2\right]\\
&=&2\,\|\phi\|\,\int_{r_1=0}^{t-s}\int_{r_2=r_1}^{t-s}\bigg(
\int_{ D}\int_{ D} q_N(s+r_1,z_1)q_N(s+r_2,z_2)p(r_1,x,z_1)p(r_2-r_1,z_1,z_2)\rho(z_1)\rho(z_2)\,dz_1\,dz_2\\
&&  -2\,\int_{ D}\int_{ \partial D} q_N(s+r_1,z_1)q(s+r_2,z_2)p(r_1,x,z_1)p(r_2-r_1,z_1,z_2)\rho(z_1)\rho(z_2)\,dz_1\,d\sigma(z_2)\\
&&  + \int_{\partial D}\int_{\partial D} q(s+r_1,z_1)q(s+r_2,z_2)p(r_1,x,z_1)p(r_2-r_1,z_1,z_2)\rho(z_1)\rho(z_2)\,d\sigma(z_1)\,d\sigma(z_2)
 \bigg)\,dr_1\,dr_2,
\end{eqnarray*}
which converges to zero uniformly for $x\in \bar{D}$ by Lemma \ref{L:MinkowskiContent_D}.
\end{pf}

\begin{remark} \label{R:2.7} \rm
While the non-negativity of $q$ easily implies that $Q$ has the contraction property (\ref{E:ContractionQNQ}), we may lose this property for $\mathbf{U}$ because intuitively the killing effect induces a jump in the system and hence can increase the fluctuation.  \qed
\end{remark}

The following gradient convergence is the cornerstone in Step 5 of the proof the main theorem. Its proof is based on the inequality $\D(P_tf)\leq (2e\,t)^{-1} \,\|f\|_{\rho}^2$ (see the Appendix of \cite{CF12}).

\begin{lem}\label{L:GradConvergenceQNQ}
For any $0\leq s\leq t$ and $\phi\in C(\bar{D})$, we have
\begin{equation}
\lim_{N\to\infty}\,\D\Big(Q^N_{(s,t)}\phi -Q_{(s,t)}\phi \Big) \,=\,0.
\end{equation}
where $\D$ is the Dirichlet form of the  $(\A,\,\rho)$-reflected diffusion defined in \eqref{E:DirichletForm_ReflectedDiffusion} and $\D  (u):=\D (u, u)$.
\end{lem}

\begin{pf}
From (\ref{E:IntegralEqt_Q}) and (\ref{E:IntegralEqt_QN}), we have
\begin{eqnarray*}
&& Q^N_{(s,t)}\phi(x)-Q_{(s,t)}\phi(x) \\
 &=& \int_0^{t-s} \int_{\partial D}p(\theta,x,y)q(s+\theta,y)Q_{(s+\theta,t)}\phi(y)\rho(y)d\sigma(y) \,-\, P_\theta (q_N(s+\theta)Q^N_{(s+\theta,t)}\phi )(x)\;d\theta\\
&=& \int_0^{t-s} P_{\theta}\left( q(s+\theta)Q_{(s+\theta,t)}\phi\,\sigma\,-\,q_N(s+\theta)Q^N_{(s+\theta,t)}\phi
\right) (x)\;d\theta\\
&=& \int_0^{t-s} P_{\theta}\left(h^{(s,t)}_N(\theta)\right)(x)\;d\theta,
\end{eqnarray*}
where $h^{(s,t)}_N(\theta)$ is the signed Borel measure $q(s+\theta,y)Q_{(s+\theta,t)}\phi(y)\,\rho(y)\sigma(dy)\,-\,q_N(s+\theta,y)Q^N_{(s+\theta,t)}\phi(y)\,\rho(y)dy$ and $P_{\theta}\mu(x):=\int_{\bar{D}}p(\theta,x,y)\,\mu(dy)$ for any measure $\mu$ on $\bar{D}$.

On the other hand, by spectral decomposition, $\D(P_tf)\leq (2e\,t)^{-1} \,\|f\|_{\rho}^2$ (see the Appendix of \cite{CF12}), where $\|\,\cdot\,\|_{\rho}$ is the $L^2(D,\,\rho(x)dx)$-norm. Hence
\begin{eqnarray}\label{E:DFTrick}
\sqrt{\D\Big(Q^N_{(s,t)}\phi(x)-Q_{(s,t)}\phi(x)\Big)} &=&
\sqrt{\D\Big(\int_0^{t-s} P_{\theta}\left(h^{(s,t)}_N(\theta)\right) (x)\;d\theta\Big)}  \notag\\
&\leq& \int_0^{t-s}\sqrt{\D\left( P_{\theta}\left( h^{(s,t)}_N(\theta)\right) \right)}\;d\theta  \notag\\
&=& \int_0^{t-s}\sqrt{\D\Big( P_{\theta/2}P_{\theta/2}\left(
h^{(s,t)}_N(\theta) \right) \Big)}\;d\theta \notag\\
&\leq&  \int_0^{t-s}\sqrt{ \frac{1}{e\,\theta}\,\Big\|P_{\theta/2}\left(
h^{(s,t)}_N(\theta)\right)\Big\|^2_{\rho}}\;d\theta.
\end{eqnarray}
We now show that the last quantity in (\ref{E:DFTrick}) converges to zero as $N\to\infty$. Note that for each $\theta\in (0,\,t-s)$, the semigroup property yields
$$\Big\|P_{\theta/2}\left( h^{(s,t)}_N(\theta)\right) \Big\|^2_{\rho}= \int_{\bar{D}}\Big(P_{\theta}h^{(s,t)}_N(\theta)\Big)(x)\,h^{(s,t)}_N(\theta)(dx) \to 0 \quad \text{as } N\to\infty$$
by Lemma \ref{L:MinkowskiContent_D} and the uniform convergence (\ref{E:UniformConvergenceQNQ}).   By the uniform bounds (\ref{E:BoundEpsilonStrip}) and (\ref{E:ContractionQNQ}), for $N$ large enough which depends only on $D$ (hence independent of $\theta$), we have $\Big\|P_{\theta}h^{(s,t)}_N(\theta)\Big\|\leq \|q\|\,\|\phi\|\,\frac{C(d,D)}{\sqrt{\theta}}$ and
\begin{eqnarray*}
\Big|P_{\theta}h^{(s,t)}_N(\theta)\Big|(\bar{D})
&=& \int_{\partial D}q(s+\theta,y)Q_{(s+\theta,t)}\phi(y)\,\rho(y)\,\sigma(dy)\,+\,\int_D q_N(s+\theta,y)Q^N_{(s+\theta,t)}\phi(y)\,\rho(y)\,dy  \\
&\leq& C(d,D)\,\|q\|\,\|\phi\|.
\end{eqnarray*}
Hence the last quantity in (\ref{E:DFTrick}) converges to zero as $N\to\infty$ by the Lebesgue dominated convergence theorem and the fact that
\begin{equation*}
\|P_{\theta/2}\mu\|_{\rho}^2 = \int_{\bar{D}}P_{\theta}\mu(x)\,\mu(dx) \leq \|P_{\theta}\mu\|\cdot|\mu|(\bar{D}),
\end{equation*}
where $|\mu|$ is the total variation measure of the signed measure $\mu$.
\end{pf}

Next, we explore the continuity in time for both $Q_{s,t}$ and $Q^N_{s,t}$.

\begin{lem}\label{L:Bound_QtMinusQs}
There exists a constant $c >0$ such that for any $0\leq s\leq t \leq T$ and $k\geq 1$,
\begin{equation*}
\sup_{r\in[0,s]}\big\|Q_{(r,t)}\phi_{k}\,-\,Q_{(r,s)}\phi_{k}\big\| \leq
c \,\|\phi_{k}\|\,\left(\lambda_{k}(t-s)+ C\,\|q\|\,(t-s)^{1/2}\right),
\end{equation*}
where $C=C(d,D,T)$ is the same constant in \eqref{E:BoundPartialD}.  Furthermore, there exists $N_0=N_0(D)$ such that for $N\geq N_0$, the above inequality holds with $\{Q^N_{s,t}\}$ in  replace of  $\{Q_{s,t}\}$.
\end{lem}

\begin{pf}
The proof will follow from a Grownwall type argument and the evolution property of the operators $\{Q_{(s,t)}\}_{s\leq t}$. By (\ref{E:IntegralEqt_Q}), for any $0\leq r\leq s\leq t$ and $k$, we have
\begin{eqnarray*}
&& \big|Q_{(r,t)}\phi_{k}(x)-Q_{(r,s)}\phi_{k}(x)\big| \\
&\leq & \big|e^{-\lambda_{k}(t-r)}\phi_{k}(x)- e^{-\lambda_{k}(s-r)}\phi_{k}(x)\big| \\
&& + \dfrac{1}{2}\Big|\int_{s-r}^{t-r}\int_{\partial D}p(\theta,x,y)q(r+\theta,y)Q_{(r+\theta,t)}\phi_{k}(y)\,\rho(y)\,d\sigma(y)\,d\theta\Big|\\
&&+ \dfrac{1}{2}\Big|\int_{0}^{s-r}\int_{\partial D}p(\theta,x,y)q(r+\theta,y)\,\big(Q_{(r+\theta,t)}\phi_{k}-Q_{(r+\theta,s)}\phi_{k}\big)(y)\,\rho(y)\,d\sigma(y)\,d\theta\Big|.
\end{eqnarray*}

Now we fix $k$, fix $0\leq s\leq t $ and define $f(r) := \big\|Q_{(r,t)}\phi_{k}-Q_{(r,s)}\phi_{k}\big\| \text{  for }r\in[0,s]$. Then the above estimate, together with (\ref{E:BoundPartialD}) and (\ref{E:ContractionQNQ}), implies that
\begin{equation}\label{E:Bound_QtMinusQs}
f(r) \leq A + B \int_0^{s-r}\dfrac{f(r+\theta)}{\sqrt{\theta}}\,d\theta \quad \text{for } r\in[0,s],
\end{equation}
where $A=\lambda_{k}\|\phi_{k}\|(t-s)+ \|q\|\,C(d,D,T)\|\phi_{k}\|(t-s)^{1/2}$ and $B=\frac{1}{2}C(d,D,T)\|q\|$.

Rewriting (\ref{E:Bound_QtMinusQs}) as $f(r) \leq A + B \int_{r}^{s}\frac{f(w)}{\sqrt{w-r}}\,dw$ and keep iterating yields
\begin{eqnarray*}
f(r) &\leq&  A+ AB\int_{w_1=r}^s\frac{1}{\sqrt{w_1-r}} + AB^2\int_{w_1=r}^s\int_{w_2=w_1}^s\frac{1}{\sqrt{(w_1-r)(w_2-w_1)}} \\
&& \; +AB^3\int_{w_1=r}^s\int_{w_2=w_1}^s\int_{w_3=w_2}^s\frac{1}{\sqrt{(w_1-r)(w_2-w_1)(w_3-w_2)}}+ \cdots  \\
&=& A\,\sum_{k=0}^{\infty} B^k\,a_k\,(s-r)^{k/2}, \text{ where } a_k=\frac{\pi^{k/2}}{\Gamma((k+2)/2)}\text{ by Lemma }\ref{L:Gamma_Half_k}\text{ in Appendix}\\
&\leq& \frac{c }{2}\,A\,\sum_{k=0}^{\infty} B^k\,(s-r)^{k/2} \quad \text{for some absolute constant }c >0\\
&\leq& c \,A \qquad \text{if}\quad |B\sqrt{s-r}|\leq 1/2
\end{eqnarray*}

Note that when $B>0$, $|B\sqrt{s-r}|\leq 1/2$ holds if and only if
$s-\frac1{4B^2}\leq s \leq s+\frac1{4B^2}$.
(The case $B=0$ is trivial since then $q=0$.)
When $0\leq r<s-1/(4B^2)$, by the evolution property and the contraction property (\ref{E:IntegralEqt_Q}), we have
\begin{eqnarray*}
\big\|Q_{(r,t)}\phi_{k} \,-\, Q_{(r,s)}\phi_{k}\big\|
&=& \big\|Q_{(r,\,s-1/(4B^2))}\left(Q_{(s-1/(4B^2),\,t)}\phi_{k} \,-\, Q_{(s-1/(4B^2),\,s)}\phi_{k}\right)\big\|\\
&\leq&  \big\| Q_{(s-1/(4B^2),\,t)}\phi_{k}\,-\,Q_{(s-1/(4B^2),\,s)}\phi_{k} \big\| \leq c \,A
\end{eqnarray*}

The above arguments clearly hold with $\{Q^N_{s,t}\}$ in  replace of  $\{Q_{s,t}\}$, if we use \eqref{E:BoundEpsilonStrip} instead of \eqref{E:BoundPartialD}. This completes the proof of the lemma.
\end{pf}

The next lemma is a key estimate that we need to establish Theorem \ref{T:ConvergenceStocInt}. Recall from (\ref{Def:D_r(phipsi)}) that
\begin{equation}\label{Def:D_r(phipsi_1)}
\D^{(q)}_r(\phi,\psi):= \<\textbf{a}\nabla \phi \cdot \nabla \psi,\,u(s)\>_{\rho} + \int_{\partial D}\phi\,\psi \,u(s)\,q(s)\,\rho\,d\sigma,\;\D^{(q)}_r(\phi):=\D^{(q)}_r(\phi,\phi).
\end{equation}
In view of \eqref{E:QuadVar_MtgMphi}, we also define
\begin{equation}\label{Def:D_r(phipsi_N)}
\D^{(q),N}_s(\phi,\psi):= \< \textbf{a}\nabla\phi\cdot\nabla\psi+q_N(s)\phi\psi,\, \X^N_s\>,\;\D^{(q),N}_s(\phi) := \D^{(q),N}_s(\phi,\phi).
\end{equation}

\begin{lem}\label{L:BoundVariance_stphi}
For all integers $k \geq 1$ and $0\leq s\leq t\leq T$, we have
\begin{eqnarray}
\int_s^t\D^{(q)}_r(Q_{(r,t)}\phi_{k})\,dr &\leq& C\,\|u_0\|\,(1\vee \|q\|)^2\,(\lambda_{k}+ \|\phi_{k}\|^2)\,(t-s),\label{E:BoundVariance_stphi_4}\\
\int_0^s\D^{(q)}_r(Q_{(r,t)}\phi_{k}-Q_{(r,s)}\phi_{k})\,dr
&\leq& C\,\|u_0\|\,(1\vee \|q\|)^4\,(\lambda_{k}^2+ \|\phi_{k}\|^2 + \|\phi_{k}\|^2\lambda_{k}^2)\,(t-s),\label{E:BoundVariance_stphi_7}
\end{eqnarray}
where $C=C(d,D,T)>0$ is a constant. Moreover, these two inequalities remain valid if we replace $Q_{r,t}$ by $Q^N_{r,t}$ and $\D^{(q)}_r$ by $\D^{(q),N}_r$ at the same time.
\end{lem}

\begin{pf}
For the first inequality, note that
\begin{equation}\label{E:BoundVariance_stphi_1}
0\leq \D^{(q)}_r(Q_{(r,t)}\phi_{k}) \leq
\|u_0\|\,\Big( \D(  Q_{(r,t)}\phi_{k})+ \sigma(\partial D)\,\|q\|\,\|\rho\|\,\|\phi_{k}^2\|\Big).
\end{equation}
Moreover, by the integral equation (\ref{E:IntegralEqt_Q}), we have
\begin{eqnarray}
\D(  Q_{(r,t)}\phi_{k}) &\leq & 2\,\D(P_{t-r}\phi_{k}) + 2\,\D\Big(\frac{1}{2}\int_{0}^{t-r}P_{\theta}\Big[H^{(r,t)}(\theta)\Big](x)\;d\theta\Big)\notag\\
&=& 2\lambda_{k}\,e^{-2(t-r)\lambda_{k}} + \D\Big(\int_{0}^{t-r}P_{\theta}\Big[H^{(r,t)}(\theta)\Big]\;d\theta\Big), \label{E:BoundVariance_stphi_2}
\end{eqnarray}
where $H^{(r,t)}(\theta)$ is the signed Borel measure $q(r+\theta,y)Q_{(r+\theta,t)}\phi_{k}(y)\,\rho(y)\,\sigma(dy)$ and $P_{\theta}\mu(x):=\int_{\bar{D}}p(t,x,y)\,\mu(dy)$ for any measure $\mu$ on $\bar{D}$.

By the same argument as that in the proof of Lemma \ref{L:GradConvergenceQNQ},
 we have
\begin{eqnarray}\label{E:BoundVariance_stphi_3}
\D\Big(\int_{0}^{t-r}P_{\theta}\Big[H^{(r,t)}(\theta)\Big]\;d\theta\Big)
&\leq& \bigg(\, \int_0^{t-r}\sqrt{\frac{1}{e\,\theta}\,\Big\|P_{\theta/2}H^{(r,t)}(\theta)\Big\|_{\rho}^2 }\;d\theta \,\bigg)^2 \notag\\
&\leq& \bigg(\, \int_0^{t-r}\sqrt{\frac{1}{e\,\theta}\,\Big\|P_{\theta}H^{(r,t)}(\theta)\Big\|_{\infty}\,|H^{(r,t)}(\theta)|(\bar{D}) }\;d\theta \,\bigg)^2 \notag\\
&\leq& \bigg(\, \int_0^{t-r} C(d,D,T)\,\|q\|\,\|\phi_{k}\|\, \theta^{-3/4} \;d\theta \,\bigg)^2 \notag\\
&\leq& C(d,D,T)\,\|q\|^2\,\|\phi_{k}\|^2\,(t-r)^{1/2}.
\end{eqnarray}
Now we put (\ref{E:BoundVariance_stphi_3}) into (\ref{E:BoundVariance_stphi_2}) and then put the result into (\ref{E:BoundVariance_stphi_1}) to obtain
\begin{equation}\label{E:BoundVariance_stphi_3.5}
\D^{(q)}_r(Q_{(r,t)}\phi_{k}) \leq \|u_0\|\,\left(\,2\lambda_{k}e^{-2(t-r)\lambda_{k}} + C(d,D,T)\,\Big(\|q\|^2\,\|\phi_{k}\|^2\,(t-r)^{1/2} +\|q\|\,\|\phi_{k}^2\|\Big)\,\right).
\end{equation}
By integration, we obtain
\begin{equation*}
\int_s^t\D^{(q)}_r(Q_{(r,t)}\phi_{k})\,dr \leq C(d,D,T)\,\|u_0\|\,\left( \|\phi_{k}\|^2\,\|q\|^2\,(t-s)^{3/2}+(\lambda_{k}+\|\phi_{k}^2\|\,\|q\|)(t-s) \right)
\end{equation*}
which implies \eqref{E:BoundVariance_stphi_4}.

The second inequality in the lemma can be dealt with in a similar way. More precisely, we have as in (\ref{E:BoundVariance_stphi_1}),
\begin{eqnarray}\label{E:BoundVariance_stphi_5}
 0 &\leq& \D^{(q)}_r(Q_{(r,t)}\phi_{k}-Q_{(r,s)}\phi_{k}) \nonumber \\
&\leq& \|u_0\|\,\left( \D(Q_{(r,t)}\phi_{k}-Q_{(r,s)}\phi_{k})
+ \sigma(\partial D)\,\|q\|\,\|\rho\|\,\big\|Q_{(r,t)}\phi_{k}-Q_{(r,s)}\phi_{k}\big\|^2\right)
\end{eqnarray}
and
\begin{eqnarray}\label{E:BoundVariance_stphi_6}
&&  \D(Q_{(r,t)}\phi_{k}-Q_{(r,s)}\phi_{k}) \notag \\
&\leq&  2 \,\left(e^{-(t-r)\lambda_{k}}-e^{-(s-r)\lambda_{k}}\right)^2\,\D(\phi_{k}) \notag\\
&&+ 2\,\D\Big(\int_{s-r}^{t-r}\int_{\partial D} p(\theta,x,y)q(r+\theta,y)Q_{(r+\theta,t)}\phi_{k}(y)\,d\sigma(y)\,d\theta\Big) \notag\\
&&+ 2\,\D\Big(\int_{0}^{s-r}\int_{\partial D} p(\theta,x,y)q(r+\theta,y)\, \big(Q_{(r+\theta,t)}\phi_{k}-Q_{(r+\theta,s)}\phi_{k}\big)(y)\,d\sigma(y)\,d\theta \Big) \notag\\
&\leq &  2 \left(e^{-(t-r)\lambda_{k}}-e^{-(s-r)\lambda_{k}}\right)^2\,\lambda_{k} \notag\\
&&+ C(d,D,T)\,\|q\|^2\,\|\phi_{k}\|^2\,(t-s)\,\left(\frac{1}{\sqrt{s-r}}-\frac{1}{\sqrt{t-r}}\right) \notag\\
&&+ C(d,D,T)\,\|q\|^2\,\left(\sup_{r\in[0,s-\theta]}\big\|Q_{(r+\theta,t)}\phi_{k}-Q_{(r+\theta,s)}\phi_{k}\big\|\right)^2\,(s-r)^{1/2} \notag\\
&\leq &  2 \left(e^{-(t-r)\lambda_{k}}-e^{-(s-r)\lambda_{k}}\right)^2\,\lambda_{k} \notag\\
&&+ C(d,D,T)\,\|q\|^2\,\|\phi_{k}\|^2\,(t-s)\,\bigg[\left(\frac{1}{\sqrt{s-r}}-\frac{1}{\sqrt{t-r}}\right)+(\lambda_k^2+\|q\|^2)(s-r)^{1/2}\bigg].
\end{eqnarray}
In the second last inequality, we have applied the same argument that we used to obtain (\ref{E:BoundVariance_stphi_3}). In the last inequality, we have used Lemma \ref{L:Bound_QtMinusQs}.

Now we put (\ref{E:BoundVariance_stphi_6}) into (\ref{E:BoundVariance_stphi_5}) and then apply Lemma \ref{L:Bound_QtMinusQs} to obtain
\begin{eqnarray*}
&& \int_0^s\D^{(q)}_r(Q_{(r,t)}\phi_{k}-Q_{(r,s)}\phi_{k})\,dr \notag \\
&\leq& \|u_0\|\,\bigg(\;
(1-e^{-(t-s)\lambda_{k}})^2(1-e^{-2s\lambda_{k}})
  + C(d,D,T)\|q\|^2\,\|\phi_{k}\|^2\,(t-s)^{3/2} \notag\\
&& \qquad\quad + C(d,D,T)\|q\|^2\,\big(\sup_{r\in[0,s]}\big\|Q_{(r+\theta,t)}\phi_{k}-Q_{(r+\theta,s)}\phi_{k}\big\|\big)^2\,s^{3/2} \notag\\
&& \qquad\quad + C(d,D,T)\|q\|\,\big(\sup_{r\in[0,s]}\big\|Q_{(r+\theta,t)}\phi_{k}-Q_{(r+\theta,s)}\phi_{k}\big\|\big)^2\,s\;
\bigg) \notag\\
&\leq&  \|u_0\|\,\lambda_{k}^2 (t-s)^2(1 \wedge 2s\lambda_{k}) \notag\\
&& + C(d,D,T)\,\|u_0\|\,\|\phi_{k}\|^2\,\Big(
\|q\|^2\,(t-s)^{3/2}+ (\|q\|^2+1)\,\big(\lambda_{k}^2 (t-s)^2+ \|q\|^2(t-s)\big)\Big) \notag\\
&\leq&  C(d,D,T)\,\|u_0\|\,\Big(\,\lambda_{k}^2 (t-s)^2+ \|\phi_{k}\|^2\,\|q\|^2\lambda_{k}^2 (t-s)^2 + \|\phi_{k}\|^2\,(\|q\|^2+\|q\|^4)(t-s)\,\Big).
\end{eqnarray*}
This implies \eqref{E:BoundVariance_stphi_7}.

Using \eqref{E:BoundEpsilonStrip} instead of \eqref{E:BoundPartialD}, we see that the above arguments remain valid if we replace $Q_{r,t}$ by $Q^N_{r,t}$ and $\D^{(q)}_r$ by $\D^{(q),N}_r$. This completes the proof of the lemma.
\end{pf}

\begin{remark}\label{Rk:BoundVariance_stphi}\rm
From the proof above, there exists $N_0=N_0(D)$ such that, for $0\leq r\leq t\leq T$ and $N\geq N_0$, inequalities \eqref{E:BoundVariance_stphi_3.5} and \eqref{E:BoundVariance_stphi_6} remain valid if we replace $Q_{r,t}$ by $Q^N_{r,t}$ and $\D^{(q)}_r$ by $\D^{(q),N}_r$. \qed
\end{remark}

\subsection{Martingales}

We need the following result from \cite[Lemma 6.1]{zqCwtF14b}.
Note that it holds for every $x\in \overline D$.

\begin{lem}\label{L:KeyMtgReflectedDiffusion}
Suppose $X=\{X_t, t\geq 0; \P_x, x\in \overline D\}$ is an $(\A,\,\rho)$-reflected diffusion in a bounded Lipschitz domain $D$ and $f$ is in the domain of the Feller generator $Dom^{Feller}(\A)$. Then we have
$$
M(t) := f(X_t)-f(X_0)-\int_0^t \A f(X_s)\,ds
$$
is a continuous $\F^{X}_t$-martingale with quadratic variation $\<M\>_t = \int_0^t \textbf{a}\nabla f \cdot \nabla f (X_s)\,ds$
under $\P^{x}$ for any $x\in \bar{D}$. Moreover, if $X_1$ and $X_2$ are independent $(\A,\,\rho)$-reflected diffusion in $D$ and $M_{i}$ is the above $M$ with $X$ replaced by $X_i$, then the cross variation $\<M_1,\,M_2\>_t=0$.
\end{lem}

From Lemma \ref{L:KeyMtgReflectedDiffusion}, we obtain the following key martingales that we need for the study of $\X^N$. The proof is the same as that for \cite[Corollary 6.4]{zqCwtF14b} so it is omitted here.

\begin{lem}\label{L:KeyMtgRobinModel}
Fix any positive integer $N$. For any $\phi\in Dom^{Feller}(\A)$, we have under $\P^{\mu}$ for any $\mu\in E_N$,
\begin{eqnarray}
M^{\phi}_t &:=& \<\phi,\X^N_t\>-\<\phi,\X^N_0\>-\int_0^t \<\A \phi-q_N(s)\,\phi,\,\X^N_s\>\,ds \quad\text{and } \label{E:MtgMphi}\\
N^{\phi}_t &:=& \<\phi,\X^N_t\>^2- \<\phi,\X^N_0\>^2 -\int_0^t
\frac{1}{N}\<\textbf{a}\nabla \phi \cdot \nabla \phi,\,\X^{N}_s\> + 2\<\phi,\X_s\>\,\<\A\phi,\X^N_s\> \notag \\
&& \qquad\qquad\qquad\qquad\qquad\quad -2\<q_N\phi,\X^N_s\>\<\phi,\X^N_s\> +\frac{1}{N}\<q_N\phi^2,\X^N_s\>\;ds \label{E:MtgNphi}
\end{eqnarray}
are $\F^{\X^N}_t$-martingales under $\P^{\mu}$ for any $\mu\in E_N$. Moreover, $M^{\phi}_t$ has predictable quadratic variation
\begin{equation}\label{E:QuadVar_MtgMphi}
\<M^{\phi}\>_t = \frac{1}{N}\int_0^t\< \textbf{a}\nabla \phi\cdot  \nabla \phi +q_N(s)\phi^2,\,\X^N_s\>\,ds.
\end{equation}
\end{lem}

From (\ref{E:QuadVar_MtgMphi}), (\ref{E:BoundEpsilonStrip}) and Lemma \ref{L:MinkowskiContent_D}, we have for all $T>0$,
\begin{equation}\label{E:QuadVar_MtgMphi_2}
\E^{\mu}[(M^{\phi}_t)^2] \leq \frac{1}{N}\left(  8(\|\phi\|^2 + \|\A \phi\|^2\,t^2) + \|\phi^2\,q\|\,C(d,D,T)t^{1/2} \right)\quad\text{for }t\in[0,T].
\end{equation}

\section{Non-equilibrium fluctuations}

In this section, we present the proof of Theorem \ref{T:Convergence_RBM_killing_Y},
the main result of this paper. Throughout this section, Assumptions \ref{A:GeometricSetting}
and \ref{A:KillingPotential} (Killing potential) are in force.
The initial distributions of the particles are assumed to be i.i.d with
distribution $u_0(x) \rho (x) dx$ for some $u_0\in C(\bar D)$.

\subsection{Langevin equation}

This subsection represents Step 1 towards the proof of Theorem \ref{T:Convergence_RBM_killing_Y}
mentioned at the end of the Introduction.
Recall that $\Y^N_t$ is the fluctuation process defined by \eqref{e:1.4}.
We first answer question (i) in the introduction of this paper.

\begin{lem}\label{L:StateSpace_RBMkilling_Y}
Whenever $\alpha>d/2$, we have $\Y^N_t \in  \H_{-\alpha}$ for $t>0$ and $N\geq 1$.
\end{lem}

\begin{pf}
Since our system is an i.i.d. sequence of sub-processes $Z^{(N)}_i$ (see Remark \ref{Rk:sub-processes}), we easily obtain
\begin{equation}\label{E:StateSpace_RBMkilling_Y_1}
\E \left[\<\Y^N_t,\,\phi\>^2\right] = {\rm Var} (\phi(Z^{(N)}_i))
 \leq \E \left[ \phi (Z^N_1 (t))^2 \right]
  \leq \<P_t\phi^2,\,u_0\> \leq \|u_0\|\,\<\phi^2,\,1\>.
\end{equation}
Hence for $\alpha>d/2$ and $t\geq0$, by \eqref{e:2.5} and \eqref{E:WeylLaw},
\begin{equation}\label{E:StateSpace_RBMkilling_Y_2}
\E\left[ |\Y^N_t|^2_{-\alpha} \right]=\sum_{k} \E\left[ \<\Y^N_t,\,h^{(\alpha)}_{k}\>^2 \right]  \leq
\|u_0\|\,\sum_{k}(1+\lambda_{k})^{-\alpha}<\infty.
\end{equation}
Then $\Y^N_t\in \H_{-\alpha}$ a.s.
\end{pf}

\begin{remark}
The condition $\alpha>d/2$ in the above lemma is sharp since, in view of
 \eqref{E:WeylLaw}, $\E\left[ |\Y^N_t|^2_{0} \right]=\infty$ when $u_0=1$, $q=0$ and $\alpha\leq d/2$.  \qed
\end{remark}

Unlike $\mathbf{U}_{(t,s)}$, we can check that $\{\mathbf{U}^N_{(t,s)}\}_{t\geq s}$ is a strongly continuous evolution system on $\H_{\gamma}$ with generator $\{\mathbf{A}^{(\gamma)}+\mathbf{B}^{(N)}_t\}_{t\geq 0}$, where $\mathbf{A}^{(\gamma)} \mu(\phi)=\mu(\A\phi)$ and $\mathbf{B}^{(N)}_t\mu(\phi)=\mu(q_N(t)\phi)$; see \cite{rfChZ95}. Using the fact that $\A\phi_{\ell}=-\lambda_{k}\,\phi_{k}$, we have $\big|\mathbf{A}^{(\gamma)}\mu\big|_{\gamma}^2 = \sum_{k}(1+\lambda_{k})^{\gamma}\lambda_{k}^2\<\mu,\phi_{k}\>^2$, which is finite if and only if (by Weyl's law) $|\mu|_{\gamma+2}^2$ is finite. Hence $Dom(\mathbf{A}^{(\gamma)})=\H_{\gamma+2}$. Since $q_N$ is bounded for each fixed $N$, we have, as operators on $\H_{\gamma}$,
\begin{equation}
Dom(\mathbf{A}^{(\gamma)}+\mathbf{B}^{(N)}_t)
=Dom(\mathbf{A}^{(\gamma)})=\H_{\gamma+2} \quad\text{for all }N\geq 1.
\end{equation}
Moreover,
\begin{equation}\label{E:UN_Contraction_RBM_Killing}
\big|\mathbf{U}^N_{(t,s)}\mu\big|_{\gamma}^2 \leq e^{(t-s)\beta_N} |\mu|^2_{\gamma}\quad \text{for some }\beta_N>0.
\end{equation}

The next result says that  $\Y^N$ solves a  \textbf{stochastic evolution equation}
in $\H_{-\alpha}$
\begin{equation}\label{E:SEE_RBM_YN_killing}
dY_t = (\mathbf{A}^{(-\alpha)}+\mathbf{B}^{(N)}_t )Y_t\,dt + dM^N_t, \quad Y_0=\Y^N_0.
\end{equation}

\begin{thm}\label{T:3.3}
Suppose $\alpha >d\vee (d/2+1)$.
For large enough $N$, there exists a r.c.l.l. square-integrable $\H_{-\alpha}$-valued martingale $M^N=(M^N_t)_{t\geq 0}$ such that $\Y^N$ satisfies the following two equivalent statements:
\begin{enumerate}
\item [(i)]\textbf{(Weak solution) } For any $\phi\in \H_{-\alpha+2}$ and $t\geq s\geq 0$, we have $\P$-a.s.
    \begin{eqnarray}\label{E:weaksol_RBM_killing}
        \<\Y^{N}_t,\,\phi\>_{-\alpha}&=& \<\Y^{N}_s,\,\phi\>_{-\alpha}+ \int_s^t \<(\mathbf{A}^{(-\alpha)}+\mathbf{B}^{(N)}_r)\Y^{N}_r,\,\phi\>_{-\alpha}\,dr + \<M^N_t-M^N_s,\,\phi\>_{-\alpha}  .
    \end{eqnarray}
\item[(ii)]\textbf{(Evolution solution) }  For $t\geq s\geq 0$, we have $\P$-a.s.
    \begin{eqnarray}\label{E:EvolutionSol_RBM_killing}
        \Y^{N}_t= \mathbf{U}^N_{(t,s)}\Y^{N}_s+ \int_s^t \mathbf{U}^N_{(t,r)}\,dM^N_r  \qquad \text{in }\H_{-\alpha}.
    \end{eqnarray}
\end{enumerate}
Moreover, $M^N$ has bounded jumps and, for every $\phi \in \H_\alpha$,
$M^N (\phi )$ is a real-valued square-integrable martingale with
 $M^N_t (\phi ) -M^N_{t-} (\phi)= \< \X^N_t-\X^N_{t-}, \phi \>$
 and predictable quadratic variation
    \begin{eqnarray}\label{E:weaksol_RBM_killing_2}
    \< M^N(\phi)\>_t  =    \int_0^t \<\textbf{a}\nabla\phi \cdot \nabla \phi\,+\,q_N(s)\phi^2,\;\X^{N}_s\>\;ds.
    \end{eqnarray}
\end{thm}

\begin{remark} \rm
Here $\int_0^{\cdot} \mathbf{U}^N_{(t,s)}\,dM^N_s$ is the stochastic integral of the operator-valued function $s\mapsto \mathbf{U}^N_{(t,s)}$ with respect to  $M^N$ on $[0,t]$. Its construction and its basic properties can be found in the monograph \cite{MP80} of M. Metivier and J. Pellaumail (See also the book by P. Protter \cite{peP05} for a more recent and comprehensive treatment for stochastic integration which used the same approach). Be aware that $t\mapsto \int_0^{t}\mathbf{U}^N_{(t,s)}\,dM^N_s$ is \emph{not} a martingale. However, since $M^N$ has a r.c.l.l. version and by (\ref{E:UN_Contraction_RBM_Killing}), we have
$\int_0^{t}\mathbf{U}^N_{(t,s)}\,dM^N_s$ has a r.c.l.l. version by the submartingale type inequality of Kotelenez (cf. \cite{pK82}).  \qed
\end{remark}

\begin{pf}
(i) and (ii) assert that $\Y^N$ is a \textbf{weak solution} and an \textbf{evolution solution} of (\ref{E:SEE_RBM_YN_killing}), respectively. Since $Dom(\mathbf{A}^{(-\alpha)})=\H_{-\alpha+2}$ is dense in $\H_{-\alpha}$, these two notion of solutions are equivalent by variation of constant (see Section 2.1.2 of \cite{GT95}). So it suffices to prove (i).

By Lemma \ref{L:KeyMtgRobinModel}, for every $\phi \in Dom^{Feller}(\A)$,
    \begin{equation}\label{E:weaksol_RBM_killing_0}
        \<\Y^N_t,\phi\>=\<\Y^N_0,\phi\>+ \int_0^t\<\Y^N_s,\,\A\phi- q_N(s)\phi\>\,ds + {M}^{N}_t(\phi),
    \end{equation}
where $ {M}^N_t(\phi)$ is a real valued $\F^{\X^{N}}_t$-martingale with quadratic variation given by (\ref{E:weaksol_RBM_killing_2})

Note that in view of \eqref{E:EigenfcnUpperBound}, each eigenfunction $\phi_k$ is bounded and continuous on $\overline D$ and hence is in the Feller generator of $\A$.
By Doob's inequality, \eqref{e:2.5}, \eqref{E:weaksol_RBM_killing_2}
 and the fact that $\E\<\phi,\X^N_s\> \leq \<P_s|\phi|,\,u_0\>$, we have
\begin{eqnarray*}
&& \sum_{k}\E\left[\sup_{[0,T]}\left(\,{M}^N_t(h^{(\alpha)}_{k})\,\right)^2 \right]\\
&\leq& C(T)\,\sum_{k}\int_0^T \E \left[ \< \textbf{a}\nabla h^{(\alpha)}_{k}\cdot \nabla h^{(\alpha)}_{k} \,+\, q_N(s)\,\big(h^{(\alpha)}_{k}\big)^2,\;\X^N_s \> \right]\;ds \\
&=& C(T)\,\sum_{k} (1+\lambda_{k})^{-\alpha} \int_0^T \<\textbf{a}\nabla \phi_{k}\cdot \nabla \phi_{k}+ q_N(s)\phi^2_{k},\;P_su_0\>_{\rho}\;ds.
\end{eqnarray*}
Recall that $\int_{\partial D} \phi_{k} (x)^2\,\sigma (dx)\leq C(d,D)(\lambda_{k}+1)$ by (\ref{E:EigenfcnUpperBound}). Hence
\begin{eqnarray}\label{L:limsup_mtg}
&& \limsup_{N\to\infty}\sum_{k}\E\left[\sup_{[0,T]}\left(\, {M}^N_t(h^{(\alpha)}_{k})\,\right)^2 \right] \nonumber \\
&\leq& C(T)\,\|u_0\|T\,\sum_{k}(1+\lambda_{k})^{-\alpha}  \Big( \D(\phi_{k})+  C(d,D)\|q\| (\lambda_{k}+1)\,\Big) \notag\\
&=& C(T)\,\|u_0\|T\,\sum_{k}(1+\lambda_{k})^{-\alpha}  \Big(\lambda_k + C(d,D)\|q\| (\lambda_{k}+1)\,\Big) \notag\\
&\leq& C(d,D,T)\,\|u_0\|\,(1\vee \|q\|)\sum_{k} \frac{1}{(1+\lambda_{k})^{\alpha-1}}
\end{eqnarray}
which by \eqref{E:WeylLaw} is finite if and only if $\alpha>d/2+1$. Hence for $\alpha>d/2+1$, there is  $N_0\geq 1$ so that for every $N\geq N_0$,
\begin{equation} \label{e:3.11}
c_N:=  \sum_{k}\E\left[\sup_{[0,T]}\left(\, {M}^N_t(h^{(\alpha)}_{k})\,
\right)^2 \right]  <\infty.
\end{equation}
For $\phi \in \H_\alpha$, $\phi = \sum_{k=1}^\infty
a_k h^{(\alpha)}_{k}$, where $a_k=\< \phi, h^{(\alpha)}_{k}\>_\alpha$.
Define ${M}^N_t(\phi)= \sum_{k=1}^\infty a_k {M}^N_t(h^{(\alpha)}_{k})$,
which is well defined in view of \eqref{e:3.11}. Moreover,
by the Doob's maximal inequality, $M^N_t(\phi)$ is the $L^2$ and
uniform limit in $t\in [0, T]$ of  $\sum_{k=1}^j a_k {M}^N_t(h^{(\alpha)}_{k})$.
Hence $M^N(\phi)$ is a real-valued r.c.l.l. square-integrable martingale with
\begin{equation}\label{e:3.12}
\E \left[ (M^N_T (\phi ))^2\right] \leq c_N \sum_{k=1}^\infty a_k^2
= c_N \| \phi\|_\alpha^2.
\end{equation}
Thus $\<M^N, \phi\>:=M^N(\phi)$ with $\phi \in \H_\alpha$ determines a r.c.l.l. square-integrable $\H_{-\alpha}$-valued martingale $M^N$. On other hand,
\begin{eqnarray*}
\sup_{t\in[0,\infty)}|M^N_t-M^N_{t-}|_{-\alpha}^2 &=&\sup_{t\in[0,\infty)}\sum_{k}(1+\lambda_{k})^{-\alpha}\,\Big(M^N_t(\phi_{k})-M^N_{t-}(\phi_{k})\Big)^2 \\
&=&\sup_{t\in[0,\infty)}\sum_{k}(1+\lambda_{k})^{-\alpha}\,N\,\Big(\X^N_t(\phi_{k})-\X^N_{t-}(\phi_{k})\Big)^2\\
&\leq& \frac{1}{N}\,\sum_{k}(1+\lambda_{k})^{-\alpha}\|\phi_{k}\|^2\\
&\leq & C/N \qquad \text{ by }
\eqref{E:WeylLaw},\;\eqref{E:EigenfcnUpperBound} \hbox{ and the assumption }\alpha>d.
\end{eqnarray*}
This in particular implies that $M^N_t$ has bounded jumps.

Finally, since $Dom(\mathbf{A}^{(-\alpha)})=\H_{-\alpha+2}$,
(\ref{E:weaksol_RBM_killing}) follows from (\ref{E:weaksol_RBM_killing_0})  provided that $\alpha> d/2+1$. This completes the proof.
\end{pf}

\subsection{Convergence of $M^N$ and tightness of $\Y^N$ }

This subsection represents Step 2 and Step 3 towards the proof of Theorem \ref{T:Convergence_RBM_killing_Y}.
By Prohorov's theorem, a sequence of $\H_{-\alpha}$-processes $\{R_N\}$ is tight in $D([0,T],\H_{-\alpha})$ provided that it satisfies the two conditions below:
  \begin{itemize}
  \item[(1)]    For all $t\in[0,T]$ and $\eps_0>0$, there exists $K>0$ such that
        \begin{equation}\label{E:Prohorov_H_MinusAlpha_1}
            \varlimsup_{N\to\infty}\P\left(|R_N(t)|^2_{-\alpha}>K\right)<\eps_0
        \end{equation}
  \item[(2)]    For all $\eps_0>0$,  as $\delta \to 0$ we have
        \begin{equation}\label{E:Prohorov_H_MinusAlpha_2}
            \varlimsup_{N\to\infty}\P\left(\sup_{\substack{|t-s|<\delta\\0\leq s,t\leq T}}\Big|R_N(t)-R_N(s)\Big|^2_{-\alpha} >\eps_0\right) \to 0
        \end{equation}
  \end{itemize}
Moreover, (\ref{E:Prohorov_H_MinusAlpha_2}) implies that any limit point has its law concentrates on $C([0,T],\H_{-\alpha})$. The following ``weak tightness criterion" can be easily checked by using (\ref{E:Prohorov_H_MinusAlpha_1}), (\ref{E:Prohorov_H_MinusAlpha_2}), the Chebyshev's inequality, the metric of $\H_{-\alpha}$.

\begin{lem}\label{L:WeakTightness_HminusAlpha}
Suppose $\{R_N; N\geq 1\}$ is a sequence of $\H_{-\alpha}$-processes for some $\alpha\in\R$ such that for any $\eps_0>0$,
\begin{equation}\label{E:WeakTightness_HminusAlpha}
  \varlimsup_{N\to\infty}\P\left(\sup_{t\in[0,T]}\,\sum_{|k|>K}\<R_N(t),\,h^{(\alpha)}_{k}\>^2\,>\eps_0\right) \to 0\quad \text{as }K\to\infty.
\end{equation}
Then the tightness of $\{R_N\}$ in $D([0,T],\H_{-\alpha})$ follows from the tightness of the one-dimensional processes $\{\<R_N,\,h^{(\alpha)}_{k}\>\}_{N\geq 1}$ (for all $k\in \mathbb{N}^d$).
\end{lem}

The following result is Step 2  towards the proof of Theorem \ref{T:Convergence_RBM_killing_Y}.

\begin{thm}\label{L:ConvergenceOfM^N_RBM_killing}
When $\alpha >d \,\vee (d/2+1)$, the square-integrable martingale $M^N$ in Theorem \ref{T:3.3} converges to $M$ in distribution in $D([0,\infty),\H_{-\alpha})$ as $N\to \infty$, where $M$ is the (unique in distribution) continuous $\H_{-\alpha}$-valued square-integrable  Gaussian martingale with independent increments and covariance functional characterized by \eqref{E:covariance_M_Killling}.
\end{thm}

\begin{pf}
We first prove the existence and uniqueness of $M$. Recall the bilinear forms $\D^{(q)}_t$ defined by (\ref{Def:D_r(phipsi)}). Fix $\alpha>d\,\vee (d/2+1)$ and define a self-adjoint  operator $A(t)$ on $\H_{-\alpha}$ by
\begin{equation}
\<A(t)\varphi^{\ast},\,\psi^{\ast}\>_{-\alpha} = \int_{0}^{t}\D^{(q)}_s(J(\varphi^{\ast}),\,J(\psi^{\ast}))\,ds,
\end{equation}
where $J:\,\H_{-\alpha}\rightarrow\H_{\alpha}$ denote the Riesz representation, i.e. for $\varphi^{\ast}\in \H_{-\alpha}$ and $\psi\in \H_{\alpha}$, we have $\<\varphi^{\ast},\,\psi\> = \<\psi,\,J(\varphi^{\ast})\>_{\alpha}$. Then $A(t)$ is a self-adjoint compact operator on the Hilbert space $\H_{-\alpha}$ of finite trace
because
$$
\sum_{k}\<A(t)h^{(-\alpha)}_{k},\,h^{(-\alpha)}_{k}\>_{-\alpha}
=\sum_{k} \int_0^t \D^{(q)}_s ( h^{(\alpha)}_{k},\,h^{(\alpha)}_{k} ) ds
<\infty
$$
 by a  calculation similar to \eqref{L:limsup_mtg}.
Moreover, $\<A(t)\varphi^{\ast},\,\varphi^{\ast}\>_{-\alpha}$ is a positive-definite quadratic functional of $\varphi^{\ast}$ for every $t$, and is continuous and increasing in $t$ for every $\varphi^{\ast}$.

Hence (cf. \cite{kI80} for a proof using Kolmogorov's extension theorem) there is a unique (in distribution) $\H_{-\alpha}$-valued Gaussian process $M$ on some probability space $(\tilde{\Omega},\,\tilde{\mathcal{F}},\,\tilde{\P})$ with independent increments, continuous sample paths, and characteristic functional
\begin{equation}
\tilde{\E}\,\exp{(i\,\<M_t,\,\varphi^{\ast}\>_{-\alpha})}= \exp{\left( -\frac{1}{2}\<A(t)\varphi^{\ast},\,\varphi^{\ast}\>_{-\alpha} \right)} .
\end{equation}

The tightness of $\{M^N\}$ and continuity of any limit are implied by Lemma \ref{L:WeakTightness_HminusAlpha} and (\ref{L:limsup_mtg}).
 Hence we only need to identify any subsequential limit. Observe that $\P$-a.s. we have
\begin{equation}\label{E:JumpSizeM^N}
\sup_{t\in[0,T]}\Big|M^N_t(\phi)-M^N_{t-}(\phi)\Big| =\sup_{t\in[0,T]}\sqrt{N}\Big|\X^N_t(\phi) - \X^N_{t-}(\phi)\Big| \leq  \frac{1}{\sqrt{N}}\,\|\phi\| \to 0
\end{equation}
and that by Theorem \ref{T:HydrodynamicLimit_RBM_Killing}, the quadratic variation (\ref{E:weaksol_RBM_killing_2}) of $M^N_t(\phi)$ converges to the deterministic quantity (\ref{E:covariance_M_Killling}) in probability for any $t\geq 0$.
 These two observations imply, by a standard functional central limit theorem for semi-martingales (see, e.g., \cite{LS80}), that $\{M^N(\phi)\}$ converges to $M(\phi)$ in distribution in $D([0,T],\R)$ for any $\phi\in Dom^{Feller}(\A)$. Finally, since $\H_{\alpha}$ has a countable dense subset in $Dom^{Feller}(\A)$ (for example, the linear span of eigenfunctions), and since any subsequential limit of $M^N$ is continuous in $t$, we know that the subsequential limit is indeed $M$. The proof is now complete.
\end{pf}

Here is Step 3  towards the proof of Theorem \ref{T:Convergence_RBM_killing_Y}.

\begin{thm}\label{T:3.7}
The sequence of processes $\{\Y^N\}$ is tight in $D([0,T],\,\H_{-\alpha})$ whenever $\alpha>d\vee (d/2+2)$. Moreover, any subsequential limit has a continuous version.
\end{thm}

\begin{pf}
We first verify (\ref{E:WeakTightness_HminusAlpha}) for $Y^N$. By (\ref{E:weaksol_RBM_killing_0}), we have
\begin{eqnarray*}
\E\Big[\sup_{[0,T]}\<\Y^N_t,\,\phi\>^2 \Big]
\leq C(T)\,\E\,\Big[\<\Y^N_0,\phi\>^2 + \int_0^T\<\Y^N_s,\A \phi\>^2ds +
\Big(\int_0^T\<\Y^N_s,\,q_N(s)\phi\>\,ds\Big)^2+\sup_{[0,T]}\left(M^N_t(\phi)\right)^2 \Big].
\end{eqnarray*}
Observe that we have treated the second term and the third term (which involve $q_N$) in the right hand side in a different way. This is because
$\int_0^T\E\<\Y^N_s,\,q_N(s)\>^2\,ds$ tends to infinity when $q$ and $u_0$ are strictly positive. The first two terms in the right hand side can be estimated using the fact
$\E [ \<\Y^N_s,\phi\>^2 ] \leq \|u_0\|\,\<\phi^2,1\>$ proved in (\ref{E:StateSpace_RBMkilling_Y_1}). The martingale term can be estimated as in (\ref{L:limsup_mtg}). For the third term which involve $q_N$, using the fact that $(\int_s^tf(r)\,dr)^2= 2\int_s^t\int_u^t f(u)f(v)\,dv\,du$ and (\ref{E:BoundEpsilonStrip}), we can check that
\begin{equation}\label{E:Tightness_RBM_killing}
\E \left[ \left(\int_s^t\<\Y^N_r,\,q_N(r)\phi\>\,dr\right)^2 \right] \leq C(D,T)\,\|\phi\|^2\,(t-s)^{3/2} \quad \text{for  }N\geq N_0(D).
\end{equation}
Combining the above calculations, we have
\begin{equation*}
\varlimsup_{N\to\infty}\sum_{k>K} \E\left[\sup_{[0,T]}\<\Y^N_t,\,h^{(\alpha)}_{k}\>^2 \right]\leq C(D,T)\,\|u_0\|\,\sum_{k>K}\,\frac{1+ \lambda_{k}^2+\|\phi_{k}^2\|+\lambda_{k}}{(1+\lambda_{k})^{\alpha}}
\end{equation*}
which,  by   (\ref{E:EigenfcnUpperBound}) and Weyl's law \eqref{E:WeylLaw}, tends to 0 as $K\to\infty$, provided that $\alpha>d\vee (d/2+2)$.
We conclude by Chebyshev's inequality that (\ref{E:WeakTightness_HminusAlpha}) for $\Y^N$ (in place of $R_N$) is satisfied if $\alpha>d\vee (d/2+2)$.

By Lemma \ref{L:WeakTightness_HminusAlpha}, it remains to verify that the one-dimensional processes $\left\{\<\Y^N,\,\phi_{k}\>; N\geq 1\right\}$ (for all $k\in \mathbb{N}$) are tight.
Since $\E [ \<\Y^N_t,\,\phi\>^2 ] \leq \|u_0\|\,\<\phi^2,\,1\>$ by (\ref{E:StateSpace_RBMkilling_Y_1}), it is enough to show that
\begin{equation}\label{E:Tightness_RBM_killing_2}
  \varlimsup_{N\to\infty}\P\left(\sup_{\substack{|t-s|<\delta \\ 0\leq s,t\leq T}}\,\big|\< \Y^N_t,\,\phi_{k}\> - \< \Y^N_s,\,\phi_{k}\>\big|\,>\eps_0\right) \to 0\quad \text{as }\delta \to 0
\end{equation}
for any $k\in \mathbb{N}$. Note that (\ref{E:Tightness_RBM_killing_2}) together with (\ref{E:WeakTightness_HminusAlpha}) for $\Y^N$ imply that any subsequential limit of $\{\Y^N; N\geq 1\}$ has a continuous version. Since $\A\phi_{k}$ is uniformly bounded and $\wh {M}^{N}_t(\phi_{k})$
defined by (\ref{E:weaksol_RBM_killing_0})  converge in $D([0,T],\R)$ as $N\to\infty$  by Theorem \ref{L:ConvergenceOfM^N_RBM_killing}, it remains to show that
\begin{equation}\label{E:Tightness_RBM_killing_3}
  \varlimsup_{N\to\infty}\P\left(\sup_{\substack{|t-s|<\delta \\ 0\leq s,t\leq T}}\,\Big|\int_s^t\<\Y^N_r,\,q_N(r)\phi_{k}\>\,dr \Big|\,>\eps_0\right) \to 0\quad \text{as }\delta \to 0.
\end{equation}
For this,  note that even though $\int_0^T\E\<q_N(s),\,\Y^N_s\>^2\,ds$ tends to infinity when $q$ and $u_0$ are strictly positive, we have
\begin{equation}\label{E:Tightness_RBM_killing}
\E\left[\,\left(\int_s^t\<\Y^N_r,\,q_N(r)\phi\>\,dr\right)^2\,\right] \leq C(T,D)\,(t-s)^{3/2} \quad \text{for  }N\geq N_0(D).
\end{equation}
This can be checked by using the fact that $(\int_s^tf(r)\,dr)^2= 2\int_s^t\int_u^t f(u)f(v)\,dvdu$. Hence we have (\ref{E:Tightness_RBM_killing_3}). See, for example, Problem 4.11 in Chapter 2 of \cite{KS91}.
\end{pf}

\subsection{Convergence of transportation part}

The goal of this subsection is to prove the following result, which is
Step 4 towards the proof of Theorem \ref{T:Convergence_RBM_killing_Y}.

\begin{thm}\label{T:3.8}
For $\alpha> d+2$,  as $N\to\infty$
\begin{equation}
\mathbf{U}^N_{(t,0)}\Y^{N}_0 \toL \mathbf{U}_{(t,0)}\Y_0 \quad \text{in  }C([0,T],\H_{-\alpha}).
\end{equation}
Moreover, $\mathbf{U}_{(t,0)}\Y_0$ has a version in $C^{\gamma}([0,T],\,\H_{-\alpha})$ for any $\gamma\in(0,1/2)$.
\end{thm}

\begin{pf}
\textbf{(i) Continuity of the limit. }We first prove that $\mathbf{U}_{(\cdot,0)}\Y_0$ has a version in $C^{\gamma}([0,T],\,\H_{-\alpha})$ for any $\gamma\in(0,1/2)$. Precisely, we will show that for $\alpha>d+2$ and $n\in\mathbb{N}$,
\begin{equation}\label{E:CtsOfUY0}
\E \left[ \big|\mathbf{U}_{(t,0)}\Y_0-\mathbf{U}_{(s,0)}\Y_0 \big|_{-\alpha}^{2n} \right] \leq  C\,\|u_0\|^n\,\left((t-s)^{2n}+ \|q\|^{2n}(t-s)^n\right)
\quad\text{whenever }0\leq s\leq t\leq T,
\end{equation}
where $C=C(n, d, D, T,\alpha)>0$ is a constant independent of $s$ and $t$. By Kolmogorov continuity criteria, (\ref{E:CtsOfUY0}) implies the desired H\"older continuity.

From Lemma \ref{L:Bound_QtMinusQs}, we have
\begin{eqnarray*}
\E \left[ \<\mathbf{U}_{(t,0)}\Y_0-\mathbf{U}_{(s,0)}\Y_0,\,\phi_{k}\>^2\right]
 &\leq&
 \| u_0 \rho \|
 \,\<\left(Q_{0, t}\phi_{k}-Q_{(0,s)}\phi_{k}\right)^2,\,1\> \\
&\leq &
\|u_0\| \,C(d,D,T)\,\|\phi_{k}\|^2\,\left(\lambda_{k}^2(t-s)^2+ \|q\|^2\,(t-s)\right).
\end{eqnarray*}
Using the Gaussian property of $\<\mathbf{U}_{(t,0)}\Y_0-\mathbf{U}_{(s,0)}\Y_0,\,\phi\>$, the above inequality and the simple fact $(a+b)^n\leq 2^n(a^n+b^n)$, we have
\begin{eqnarray*}
&& \E\left[ \<\,\mathbf{U}_{(t,0)}\Y_0-\mathbf{U}_{(s,0)}\Y_0,\;\phi_{k}\,\>^{2n}
\right] = (2n-1)!!\,\left(\E \left[ \<\mathbf{U}_{(t,0)}\Y_0-\mathbf{U}_{(s,0)}\Y_0,\,\phi\>^2\right] \right)^n \\
&\leq & (2n-1)!!\,2^n\,C^n(d,D,T)\,\|u_0\|^n\,\|\phi_{k}\|^{2n}\,\left(\lambda_{k}^{2n}(t-s)^{2n}+ \|q\|^{2n}\,(t-s)^n\right).
\end{eqnarray*}
Therefore, using H\"older inequality $(\sum_{i} a_ib_i)^n \leq (\sum_i a_i^{n/(n-1)})^{n-1}(\sum_i b_i^n)$ for non-negative numbers $a_i$ and $b_i$, we have for any $\beta\in (0,\alpha]$,
\begin{eqnarray*}
&& \E \left[ \big|\mathbf{U}_t\Y_0-\mathbf{U}_s\Y_0 \big|_{-\alpha}^{2n} \right]\\
&=& \E \left[ \Big(\sum_{k}(1+\lambda_{k})^{-\alpha}\,\<\mathbf{U}_{t}\Y_0-\mathbf{U}_{s}\Y_0,\,\phi_{k}\>^2
\Big)^n \right] \\
&\leq& \left(\sum_{k}(1+\lambda_{k})^{-\frac{\beta n}{n-1}}\right)^{n-1}
\left(\sum_{k}(1+\lambda_{k})^{-(\alpha-\beta)n}\,\E\left[
\<\mathbf{U}_{t}\Y_0-\mathbf{U}_{s}\Y_0,\,\phi_{k}\>^{2n} \right] \right) \\
&\leq& C(n,d,D,T)\,\|u_0\|^n\,\left(\sum_{k}\frac{1}{(1+\lambda_{k})^{\frac{\beta n}{n-1}}}\right)^{n-1}\\
&& \qquad\qquad \cdot\left(\,(t-s)^{2n}\sum_{k} \frac{\|\phi_{k}\|^{2n}\,\lambda_{k}^{2n}}{(1+\lambda_{k})^{(\alpha-\beta)n}} \,+\,
\|q\|^{2n}\,(t-s)^{n}\sum_{k} \frac{\|\phi_{k}\|^{2n}}{(1+\lambda_{k})^{(\alpha-\beta)n}}\right).
\end{eqnarray*}
From (\ref{E:EigenfcnUpperBound}), it follows that (\ref{E:TightnessOfUNYN0}) holds true once we choose $\beta\in \left( \frac{d(n-1)}{2n},\,\alpha-\frac{d}{2}-2- \frac{d}{2n}\right)$. This choice of $\beta$ is possible if and only if $\alpha>d+2$. Hence the proof of (\ref{E:CtsOfUY0}) is complete.

\textbf{(ii) Tightness. }Next, we show that $\{\mathbf{U}^N_{(\cdot,0)}\Y^{N}_0\}$ is tight in $C([0,T],\H_{-\alpha})$. Let $\psi=Q^N_{0,t}\phi_k-Q^N_{0,s}\phi_k$ and $\{x_i\}_{i=1}^N$ be i.i.d. with distribution $u_0(x) \rho (x) dx$.
Then
\begin{eqnarray*}
&&\E \left[ \<\mathbf{U}^N_{(t,0)}\Y^N_0-\mathbf{U}^N_{(s,0)}\Y^N_0,\,\phi_{k}\>^4
 \right]= \E \left[ \<\Y^N_0,\,\psi\>^4 \right] \\
&=& N^2\,\E\bigg[\, \bigg(\frac{1}{N}\sum_{i=1}^N\Big(\psi(x_i)-\mu_{\psi}\Big)\bigg)^4\,\bigg]\qquad\quad \text{(where }\mu_{\psi}:=\<\psi,\,u_0\>_{\rho} \text{ is the mean of each term)}\\
&=&\frac{1}{N}\,\E[(\psi(x_1)-\mu_{\psi})^4]\,+\,\Big(\,\E[(\psi(x_1)-\mu_{\psi})^2]\,\Big)^2\\
&\leq&  C(d,D,T)\,
\|u_0 \rho \|^2
\,\|\phi_{k}\|^{4}\,\Big(\lambda_{k}^{4}(t-s)^{4}+ \|q\|^{4}\,(t-s)^2\Big)\quad \text{(by Lemma }\ref{L:Bound_QtMinusQs}\text{)}.
\end{eqnarray*}
Using H\"older inequality $(\sum_{i} a_ib_i)^n \leq (\sum_i a_i^{n/(n-1)})^{n-1}(\sum_i b_i^n)$ as in step (i) above (with $n=2$ here), we obtain
\begin{equation}\label{E:TightnessOfUNYN0}
 \sup_{N>1} \E \left[ \big|\mathbf{U}^N_{(t,0)}\Y^{N}_0-\mathbf{U}^N_{(s,0)}\Y^{N}_0 \big|_{-\alpha}^{4} \right] \leq  C\,\|u_0\|^2\,\left((t-s)^{4}+ \|q\|^{4}(t-s)^2\right)
\end{equation}
whenever $0\leq s\leq t\leq T$ and $\alpha>d+2$, where $C=C(\alpha,d,D,T)>0$. Inequality  \eqref{E:TightnessOfUNYN0} implies the tightness we need, in view of the Kolmogorov-Centov tightness criteria (see \cite[Theorem 3.8.8]{EK86}).

\textbf{(iii) Convergence of finite dimensional distributions. } To finish the proof of Theorem \ref{T:3.8}, it remains to show that for any $n\in \mathbb{N}$ and $0\leq t_1\leq\cdots \leq t_n<\infty$, we have
\begin{equation}\label{E:FiniteDimConvergence_UNYN0_Halphak}
\left(\mathbf{U}^N_{(t_1,0)}\Y^{N}_0,\,\cdots,\,\mathbf{U}^N_{(t_n,0)}\Y^{N}_0\right) \toL \left(\mathbf{U}_{(t_1,0)}\Y_0,\,\cdots,\,\mathbf{U}_{(t_n,0)}\Y_0\right) \quad \text{in  }(\H_{-\alpha})^n
\end{equation}
as $N\to\infty$.

For this, it suffices to show that for any $\psi_1,\cdots,\psi_n\in \mathcal{C}\subset C(\bar{D})$,
\begin{equation}\label{E:FiniteDimConvergence_UNYN0_Rk}
\left(\<\mathbf{U}^N_{(t_1,0)}\Y^{N}_0,\,\psi_1\>,\,\cdots,\,\<\mathbf{U}^N_{(t_n,0)}\Y^{N}_0,\,\psi_n\>\right) \toL \left(\<\mathbf{U}_{(t_1,0)}\Y_0,\,\psi_1\>,\,\cdots,\,\<\mathbf{U}_{(t_n,0)}\Y_0,\,\psi_n\>\right) \quad \text{in  }\R^n,
\end{equation}
where $\mathcal{C}$ denotes the linear span of the eigenfunctions $\{\phi_{k}\}$. This is because $\mathcal{C}$ is dense in $\H_{\alpha}$ and the Borel $\sigma$-field in $(\H_{-\alpha})^n$ is generated by the finite dimensional sets.

We first prove (\ref{E:FiniteDimConvergence_UNYN0_Rk}) when $n=1$.
For notational simplicity, write $t$ and $\psi$ for $t_1$ and $\psi_1$. Note that $\left\{\<\mathbf{U}^N_{(t,0)}\Y^{N}_0,\,\psi\>\right\}$ is tight in $\R$ since
by \eqref{E:ContractionQNQ},
\begin{eqnarray*}
 \sup_N \E \left[ \<\mathbf{U}^N_{t}\Y^{N}_0,\,\psi\>^2 \right] =\sup_N \< \big(Q_{0,t}^N\psi\big)^2,\,u_0\>- \<Q_{0,t}^N\psi,\,u_0\>^2
 \leq  \|u_0\|  \,\sup_N \<  (Q_{0,t}^N\psi )^2,\,1\>  <  \infty .
\end{eqnarray*}
Suppose $Z$ is a subsequential limit of $\<\mathbf{U}^N_{(t,0)}\Y^{N}_0,\,\psi\>$. We claim  that $Z\equalL \<\mathbf{U}_{(t,0)}\Y_0,\,\psi\>$. This is due to the following two facts: $\<\Y^{N}_0,\,Q_{0, t}\psi\> \toL \<\Y_0,\,Q_{0, t}\psi\>$ (by the standard central limit theorem) and $\lim_{N\to 0}  \E\big|\<\Y^{N}_0,\,Q^N_{0, t}\psi\> -\<\Y^{N}_0,\,Q_{0, t}\psi\> \big|
= 0$, which follows from
\begin{eqnarray*}
 \E\left[ \big|\<\Y^{N}_0,\,Q^N_{0, t}\psi\> -\<\Y^{N}_0,\,Q_{0, t}\psi\> \big|^2\right]
&=&\< ( Q_{0, t}^N\phi-Q_{0, t}\psi )^2,\,u_0\>- \<Q_{0, t}^N\psi-Q_{0, t}\psi,\,u_0\>^2\\
&\leq& \|u_0\|\,\< (Q_{0, t}^N\psi-Q_{0, t}\psi )^2,\,1\> \to 0 \quad  \text{ by }(\ref{E:UniformConvergenceQNQ}).
\end{eqnarray*}
In fact,  the second fact implies that $\<\Y^{N'}_0,\,Q^{N'}_{0,t}\psi\> -\<\Y^{N'}_0,\,Q_{0, t}\psi\> \to 0$ a.s. along some subsequence $N'$, and so by the Lebesgue dominated convergence theorem,  $\E F(\<\Y^{N'}_0,\,Q^{N'}_{0,t}\psi\>)-\E F(\<\Y^{N'}_0,\,Q_{0, t}\psi\>)\to 0$ for any bounded continuous function $F$.

The proof of (\ref{E:FiniteDimConvergence_UNYN0_Rk}) for general $n\in \mathbb{N}$ is the same as that for $n=1$, using the standard multidimensional central limit theorem. So we get the desired \eqref{E:FiniteDimConvergence_UNYN0_Halphak}.

The proof of Theorem \ref{T:3.8} is now complete.
\end{pf}

\subsection{Convergence of stochastic integrals}

Our goal in this subsection is to prove the following result,
which corresponds to  Step 5 towards the proof of Theorem \ref{T:Convergence_RBM_killing_Y}.

\begin{thm}\label{T:ConvergenceStocInt}
For $\alpha> d+2$ and $T>0$, as  $N\to\infty$
\begin{equation}
\int_0^{t} \mathbf{U}^N_{(t,s)}dM^N_{s} \toL \int_0^{t} \mathbf{U}_{(t,s)}dM_{s} \quad \text{in  }D([0,T],\H_{-\alpha}) .
\end{equation}
Moreover,  $\int_0^{t} \mathbf{U}_{(t,s)}dM_{s}$ has a version in $C^{\gamma}([0,T],\,\H_{-\alpha})$ for any $\gamma\in(0,1/2)$.
\end{thm}

First, we need the following lemma which is the key for establishing finite dimensional convergence. Lemma \ref{L:BoundVariance_stphi} also plays a crucial role in the proof of Theorem \ref{T:ConvergenceStocInt}. Recall from \eqref{Def:D_r(phipsi_1)} and \eqref{Def:D_r(phipsi_N)} that
\begin{eqnarray*}
\D^{(q)}_r(\phi,\psi)&:=& \<\textbf{a}\nabla \phi \cdot \nabla \psi,\,u(s)\>_{\rho} + \int_{\partial D}\phi\,\psi \,u(s)\,q(s)\,\rho\,d\sigma,\;\D^{(q)}_r(\phi):=\D^{(q)}_r(\phi,\phi), \text{ and}\\
\D^{(q),N}_s(\phi,\psi)&:=& \< \textbf{a}\nabla\phi\cdot\nabla\psi+q_N(s)\phi\psi,\, \X^N_s\>,\;\D^{(q),N}_s(\phi) := \D^{(q),N}_s(\phi,\phi).
\end{eqnarray*}

\begin{lem}\label{L:ConvergenceDistOfIntegral}
For $0\leq a\leq b\leq T$, $i=\sqrt{-1}$ and $\phi \in C(\bar{D})$, as
$N \to \infty$, we have
\begin{eqnarray*}
&&  \E\bigg[ \exp{ \Big(i\big\langle
                \int^b_a \mathbf{U}^N_{(T,s)}dM^N_{s},\,\phi
            \big\rangle \Big)} \Big|\F^N_a  \bigg]
\quad \hbox{converges in } L^1(\P)\hbox{ to }\\
 && \hskip 0.3truein
 \exp{\Big(-\frac{1}{2} \int^b_a \D^{(q)}_s(Q_{(s,T)}\phi)ds \Big)}
= \E\bigg[ \exp{ \Big( i\big\langle \int^b_a \mathbf{U}_{(T,s)}dM_{s},\,\phi \big\rangle \Big) } \bigg] .
\end{eqnarray*}
\end{lem}

\begin{pf}
\textbf{(i) }
 Fix $T>0$ and $\phi \in \H_\alpha$. Then
$$
K_t=K^N_t:= \big\langle\int^t_0 \mathbf{U}^N_{(T,s)}dM^N_{s},\,\phi \big\rangle \quad\text{is a martingale for }t \in [0,T].
$$
Let $\triangle K_r := K_r - K_{r-}$ denote the jump of $K$ at time $r$. Then by
(\ref{E:JumpSizeM^N}) and (\ref{E:ContractionQNQ}),
\begin{equation}\label{E:JumpSizeK}
\sup_{r\in[0,T]}|\triangle K_r| \leq \sup_{s\in[0,T]}\|Q^N_{s,T}\phi\| / \sqrt{N} \leq \|\phi\| / \sqrt{N} .
\end{equation}
Moreover, by Theorem \ref{T:3.3}, the dual predictable projection $\< K\>$
of the quadratic variation $[K]$ of $K$ is
\begin{equation}\label{E:QuadVarK}
\< K\> _t = \int^t_0 \D^{(q),N}_s(Q^N_{(s,T)}\phi)ds \qquad\hbox{for }
 t \in [0,T],
\end{equation}
By a similar argument as that for $H_N(t)$ in the proof of Lemma \ref{L:HNGNtoZero} (using an inequality in Remark \ref{Rk:BoundVariance_stphi}), we have $\limsup_{N\to\infty}\E \big[ \< K\>^k_T \big]<\infty $ for every integer $k\geq 1$.
Observe that
$ J_t :=
[K]_t - \<K\>_t$ is a purely discontinuous martingale with jumps $\Delta J_t:=J_t - J_{t-}= (\Delta K_t)^2$. It follows from (\ref{E:JumpSizeK}) that
$$
\E [ J_T^2 ] = \E \left[ \sum_{0<s\leq T} (\Delta J_s)^2 \right] \leq  \| \phi \|^2 N^{-1}  \, \E \left[ \sum_{0<s\leq T} (\Delta K_s)^2 \right]
\leq  \| \phi \|^2 N^{-1}  \, \E  [ K]_T = \| \phi \|^2 N^{-1} \,\E\< K \>_T.
$$
Hence
\begin{equation}\label{e:4.30}
\E  \left[ [K]_T^2 \right] = \E \left[ ( \<K\>_T + J_T)^2 \right] \leq 2 \E \left[ \<K\>^2_T +J^2_T \right] \leq 2\E\<K\>_T+2\| \phi \|^2 N^{-1} \,\E\< K \>_T
\end{equation}
which is uniformly bounded in $N$, by Lemma \ref{L:BoundVariance_stphi}.

Let $f(r):= e^{ir}$, $g(r) := \D^{(q)}_r(Q_{r,T}\phi)$, and $g_N(r) := \D^{(q),N}_r(Q^N_{r,T}\phi)$. Fix $a \in [0,T]$,
and set  $h_N(t) := \E\Big[f(K_t-K_a)|\F^N_a\Big]$ and $h(t) := \exp{\Big(-\frac{1}{2}\int^t_ag(r)dr\Big)}$.
Note that $h(t) = 1-\frac{1}{2}\int^t_a h(r)g(r)dr$. We claim that
\begin{equation}\label{E:Equation_h_N}
h_N(t)=1-\frac{1}{2}\int^t_a h_N(r)g(r)dr + \eps_N(t) \  \text{ with } \sup_{t\in[a,T]}|\eps_N(t)| \rightarrow 0 \,\hbox{ in   } L^1(\P).
\end{equation}
By Gronwall's inequality, the above equations yield
\begin{equation*}
|h_N(t)-h(t)| \leq \bigg(\sup_{t\in[a,T]}|\eps_N(t)|\bigg)\exp{\bigg(\frac{1}{2}\int^t_ag(r)dr\bigg)}
\end{equation*}
and hence $h_N(t) \rightarrow h(t)$ in $L^1(\P)$ as $N\to \infty$.
On other hand,  since $M$ is Gaussian with independent increment, $\big\langle \int^b_a \mathbf{U}_{(c,s)}\,dM_{s},\,\phi \big\rangle$ is Gaussian  with variance $\int^b_a \D^{(q)}_s(Q_{(s,T)}\phi)ds$ (see subsection \ref{subsection:WellDefineStocInt} in the Appendix). Thus we have
$\exp{\left(-\frac{1}{2} \int^b_a \D^{(q)}_s(Q_{(s,T)}\phi)ds \right)}
= \E\left[ \exp{ \left( i\big\langle \int^b_a \mathbf{U}_{(T,s)}\,dM_{s},\,\phi \big\rangle \right) } \right]$.
This proves Lemma \ref{L:ConvergenceDistOfIntegral} once the claim \eqref{E:Equation_h_N} is verified. We now prove (\ref{E:Equation_h_N})
 in the next two steps.

\textbf{(ii) }
By It\'{o}'s formula (see, e.g.,  Theorem 36 in  \cite[Chapter II]{peP05}),
\begin{eqnarray}\label{E:ItofKt}
f(K_t) &=& 1+\int^t_{0+} f'(K_{r-})dK_r + \frac{1}{2}\int^t_{0+}f''(K_{r-})d[K]_r \notag \\
&& +\sum_{0< r\leq t} \left( f(K_r)-f(K_{r-})-f'(K_{r-}) \triangle K_r - \frac{1}{2} f''(K_{r-})(\triangle K_r)^2\right).
\end{eqnarray}

Hence for $t \in [a,T]$,
\begin{eqnarray*}
\E\left[f(K_t) \big|\F^N_a\right] &=& f(K_a)+\E\left[\frac{1}{2}\int^t_{a+}f''(K_{r-})d[K]_r \big| \F^N_a\right]\\
&&+\E \left[ \sum_{0<r\leq t} \left( f(K_r)-f(K_{r-})-f'(K_{r-}) \triangle K_r
 - \frac{1}{2} f''(K_{r-})(\triangle K_r )^2\right)\Big|\F^N_a\right]\\
&=& f(K_a)-\frac{1}{2}\E\left[ \int^t_{a+}f(K_{r-})g(r)dr\Big|\F^N_a\right]+\eps^{(1)}_N(t)+\eps^{(2)}_N(t),
\end{eqnarray*}
where
\begin{eqnarray*}
\eps^{(1)}_N(t) &:=& \frac{1}{2}\E \left[ \int^t_{a+}f(K_{r-})(g(r)-g_N(r))dr\Big|\F^N_a\right]\quad\text{and}\\
\eps^{(2)}_N(t) &:=& \E\left[\sum_{0<r\leq t} \left( f(K_r)-f(K_{r-})-f'(K_{r-}) \triangle K_r - \frac{1}{2} f''(K_{r-})(\triangle K_r)^2\right)\Big|\F^N_a\right].
\end{eqnarray*}
We have used (\ref{E:QuadVarK}) and the fact that $f''=-f$ in the last equality.

Dividing both sides by $f(K_a)$, the above calculations give
\begin{equation}\label{e:4.33}
h_N(t)=1-\frac{1}{2}\int^{t}_{a+}h_N(r)g(r)dr + \frac{\eps^{(1)}_N(t)+\eps^{(2)}_N(t)}{f(K_a)}.
\end{equation}
Since $|f|=1$ and $|e^{ia}-1-ia+a^2/2|\leq |a|^3/6$,
we have by (\ref{E:JumpSizeK})
$$
  \left|\eps^{(2)}_N(t)\right|
 \leq  \frac{1}{6}\,\E\bigg[\sum_{0<r\leq T}|\triangle K_r|^3\Big|\,\F^N_a\bigg]
 \leq  \frac{\|\phi\|}{6\sqrt{N}} \,\E\bigg[\sum_{0<r\leq T}(\triangle K_r)^2\Big|\,\F^N_a\bigg]
 \leq  \frac{\|\phi\|}{6\sqrt{N}} \,\E\left[[K]_T \Big|\,\F^N_a\right].
$$
Since $\E [K]_T = \int_0^T \E[g_N(s)]ds\rightarrow \int^T_0g(s)\,ds$, we get $\lim_{N\to \infty} \E \left[ \sup_{t\in[a,T]} | \eps^{(2)}_N(t) | \right] =0$.

For $\eps^{(1)}_N$, we let $\psi (r) := Q_{(r,T)}\phi$ and $\psi_N (r) :=Q^N_{(r,T)}\phi$ for simplification. Since $|f|=1$, triangle inequality gives
\begin{eqnarray}\label{E:ConvergenceDistOfIntegral_2}
&& 2\sup_{t\in[a,T]}\left|\eps^{(1)}_N(t)\right| \notag \\
&\leq& \E\bigg[\int^T_{a+} \big|\<\textbf{a}\nabla\psi\cdot \nabla\psi,\,u(r)\>_{\rho}-\<\textbf{a}\nabla\psi\cdot \nabla\psi,\,\X^N_r\>\big|\,dr\,\Big|\,\F^N_a\bigg] \notag\\
&&+\,\E\bigg[\int^T_{a+} \big|\<\textbf{a}\nabla\psi\cdot \nabla\psi-\textbf{a}\nabla\psi_N\cdot \nabla\psi_N,\;\X^N_r\>\big|\,dr\,\Big|\,\F^N_a\bigg] \notag\\
&&+\,\int^T_{a+} \Big|\int_{\partial D}\psi^2q(r)u(r)\,\rho\,d\sigma-\<\psi^2_Nq_N(r),\,u(r)\>_{\rho} \Big|\,dr \notag\\
&&+\,\sup_{t\in[a,T]}\,\bigg|\E\bigg[\int^t_{a+} f(K_{r-})\Big(\<\psi^2_Nq_N(r),\,u(r)\>_{\rho}-\<\psi_N^2q_N(r),\,\X^N_r\> \Big)\,dr\,\Big|\,\F^N_a\bigg]\,\bigg|.
\end{eqnarray}

The expectation of the first term on the right hand side of (\ref{E:ConvergenceDistOfIntegral_2}) tends to zero by the hydrodynamic result (Theorem \ref{T:HydrodynamicLimit_RBM_Killing}). The expectation of the second term is at most
\begin{eqnarray*}
&&\E\int^T_{a+}\big\langle\big|\textbf{a}\nabla\psi\cdot \nabla\psi-\textbf{a}\nabla\psi_N\cdot \nabla\psi_N\big|,\;\X^N_r\big\rangle\; dr\\
&\leq& \int^T_{a+}\big\langle P_r\big(\big|\textbf{a}\nabla\psi\cdot \nabla\psi-\textbf{a}\nabla\psi_N\cdot \nabla\psi_N\big|\big),\;u_0\big\rangle_{\rho} \;dr\\
&\leq& \| u_0\| \int^T_{a+} \big\langle \big|\textbf{a}\nabla\psi\cdot \nabla\psi-\textbf{a}\nabla\psi_N\cdot \nabla\psi_N\big|,\;1 \big\rangle_{\rho} \;dr\\
&=& \| u_0\| \int^T_{a+} \big\langle \big|\textbf{a}\nabla(\psi-\psi_N)\cdot \nabla(\psi+\psi_N)\big|,\;1 \big\rangle_{\rho} \;dr \quad \text{by symmetry of }\textbf{a}\\
&\leq& \| u_0\| \int^T_{a+} \sqrt{\D(\psi-\psi_N)\,\D(\psi+\psi_N)} \;dr \quad \text{by Cauchy-Schwarz inequality }.
\end{eqnarray*}
This last quantity tends to zero as $N\to\infty$ by Lemma \ref{L:GradConvergenceQNQ} and Lebesgue dominated convergence theorem.

The third term (which is deterministic) on the right hand side of (\ref{E:ConvergenceDistOfIntegral_2}) converges to zero as $N\to\infty$ by Lemma \ref{L:MinkowskiContent_D} and the uniform convergence (\ref{E:UniformConvergenceQNQ}).

\textbf{(iii) }
It remains to show that the forth (and last) term on the right hand side of (\ref{E:ConvergenceDistOfIntegral_2}) converges to zero in $L^1(\P)$. This term can be written as
\begin{equation}\label{E:ForthTerm}
\sup_{t\in[a,T]}\,\bigg|\E\bigg[\int^t_{a+} f(K_{r-})\,dH_N(r) - \int^t_{a+} f(K_{r-})\,dG_N(r)
\,\Big|\,\F^N_a\bigg]\,\bigg|,
\end{equation}
where
$$H_N(t):=\int_{a+}^t\<\psi_N^2q_N(r),\,\X^N_r\>\,dr \quad\text{and}\quad G_N(t):=\int_{a+}^t \<\psi^2_Nq_N(r),\,u(r)\>_{\rho}\,dr.$$
We have by Lemma \ref{L:HNGNtoZero}
\begin{equation}\label{E:ForthTerm_1}
\lim_{N\to\infty} \E\left[\,\Big(\sup_{t\in[0,T]}|H_N(t)-G_N(t)|\Big)^p\,\right] =0\quad \text{for every }p\geq 1.
\end{equation}
In view of \eqref{E:ItofKt} and \eqref{E:JumpSizeK}, it suffice to show (\ref{E:ForthTerm}) converges to zero in $L^1(\P)$ with $f(K_{t})$ replaced by
$\tilde{f}(K_{t}):= 1+\int^t_{0+} f'(K_{r-})dK_r + \frac{1}{2}\int^t_{0+}f''(K_{r-})d[K]_r$. Furthermore, since $H_N(t)$ and $G_N(t)$ have bounded variations, by an integration by parts (see, e.g.,  Corollary 2 in  \cite[Chapter II]{peP05}), we have
\begin{eqnarray*}
\int^t_{a+} \tilde{f}(K_{r-})\,dH_N(r) &=& \tilde{f}(K_{t})\,H_N(t) - \int^t_{a+} H_N(r) \,d\tilde{f}(K_{r}) \quad \text{and}\\
\int^t_{a+} \tilde{f}(K_{r-})\,dG_N(r) &=& \tilde{f}(K_{t})\,G_N(t)- \int^t_{a+} G_N(r) \,d\tilde{f}(K_{r}).
\end{eqnarray*}
On subtraction, it suffices to show
\begin{equation}\label{E:ForthTerm_2}
\E\Big[ \sup_{t\in[a,T]}\, \big| \tilde{f}(K_t)\,\big(H_N(t)-G_N(t)\big) \big|\,\Big|\,\F^N_a\Big] \quad \text{and}
\end{equation}
\begin{equation}\label{E:ForthTerm_3}
\sup_{t\in[a,T]}\,\bigg|\E\bigg[\int^t_{a+} \Big(H_N(r)-G_N(r)\Big)\,f(K_r)\,d[K]_r
\,\Big|\,\F^N_a\bigg]\,\bigg|
\end{equation}
both converge to zero in $L^1(\P)$.

Since $|f|=1$, $|e^{ia}-1-ia+a^2/2|\leq |a|^3/6$, we have
by (\ref{E:JumpSizeK}) and \eqref{E:ItofKt},
$$
\sup_{t\in[a,T]}|\tilde{f}(K_{t-})| \leq \sup_{t\in[a,T]}\big(|f(K_{t-})|+ |f(K_{t-})-\tilde{f}(K_{t-})| \Big) \leq 1+ \frac{\|\phi\|\,[K]_T}{6\sqrt{N}}.
$$
Hence (\ref{E:ForthTerm_2}) converges to zero in $L^1(\P)$ by \eqref{e:4.30} and \eqref{E:ForthTerm_1}.
Finally, the expectation of (\ref{E:ForthTerm_3}) is at most
$$
 \E\left[ \sup_{t\in[a,T]} \left|   \left( H_N(r)-G_N(r)\right)\right|
 \left( [K]_T -[K]_a  \right)  \right] \leq
 \left( \E\left[ \sup_{t\in[a,T]}
 \left( H_N(r)-G_N(r)\right)^2\right] \right)^{1/2}
\left( \E\left[ [K]_T^2 \right]   \right)^{1/2} ,
$$
which goes to $0$ as $N\to \infty$ by \eqref{e:4.30} and \eqref{E:ForthTerm_1}.
Hence by \eqref{E:ConvergenceDistOfIntegral_2},
$\sup_{t\in[a,T]}|\eps^{(1)}_N(t)|\to 0$ in $L^1(\P)$. We then conclude
from \eqref{e:4.33} that \eqref{E:Equation_h_N} holds.
The proof of the lemma is now complete.
\end{pf}

\medskip

We can now present the proof of Theorem \ref{T:ConvergenceStocInt}.

\medskip

\begin{pf}\emph{of Theorem} \ref{T:ConvergenceStocInt}.
For notational convenience, set $J_N(t) := \int^t_0 \mathbf{U}^N_{(t,r)}\,dM^N_{r}$ and $J(t) := \int^t_0 \mathbf{U}_{(t,r)}\,dM_{r}$.

\textbf{(i) Continuity of the limit. }
In the Appendix, we checked that $J(t)$ is a well-defined $\H_{-\alpha}$-valued Gaussian random variable. We now prove that $J(\cdot)$ has a version in $C^{\gamma}([0,T],\,\H_{-\alpha})$ for any $\gamma\in(0,1/2)$. By Kolmogorov continuity criteria, it suffices to show that for $\alpha>d+2$ and $n\in\mathbb{N}$,
\begin{equation}\label{E:CtsOfIntUdM}
\E \Big[\, | J(t) -J(s) |_{-\alpha}^{2n} \,\Big] \leq  C\,(t-s)^{n}\quad\text{whenever }0\leq s\leq t\leq T,
\end{equation}
where $C=C(n, d, D, T,\alpha)\, \|u_0\|^n(1\vee \|q\|)^{4n}>0$ is a constant.

Note that  for $\phi\in C(\bar{D})$,
$$\<\,J(t) -J(s),\;\phi\,\> =
\Big\langle \int_s^t\mathbf{U}_{(t,r)}dM_{r},\;\phi \Big\rangle + \Big\langle \int_0^s\mathbf{U}_{(t,r)}-\mathbf{U}_{(s,r)}\,dM_{r},\;\phi \Big\rangle ,
$$
which, as the sum of two independent centered Gaussian variable,  is a centered Gaussian random variable with variance $$V_s^t(\phi_{k}):=\int_s^t\D^{(q)}_r(Q_{(r,t)}\phi_{k})\,dr+ \int_0^s\D^{(q)}_r(Q_{(r,t)}\phi_{k}-Q_{(r,s)}\phi_{k})\,dr.$$
By Lemma \ref{L:BoundVariance_stphi}, we have
\begin{eqnarray*}
&& \E\left[ \<\,J(t) -J(s)\,,\;\phi_{k}\,\>^{2n} \right]
= (2n-1)!!\,\left(V_s^t(\phi_{k})\right)^n\\
&\leq& C(n, d, D, T) \,\|u_0\|^n(1\vee \|q\|)^{4n}\,(1\vee \lambda_{k})^{2n}\,\|\phi_{k}\|^{2n}\,(t-s)^{n}
\end{eqnarray*}
for any $0\leq s\leq t \leq T$ and $k\in \mathbb{N}$. Applying H\"older inequality $(\sum_{i} a_ib_i)^n \leq (\sum_i a_i^{n/(n-1)})^{n-1}(\sum_i b_i^n)$, we have for any $\beta\in (0,\alpha]$,
\begin{eqnarray*}
&& \E \left[ |J(t) -J(s)|_{-\alpha}^{2n} \right] \\
&=& \E \left[ \left(\sum_{k}(1+\lambda_{k})^{-\alpha}\, \<\,J(t) -J(s)\,,\;\phi_{k}\,\>^{2}\right)^n \right] \\
&\leq& \left(\sum_{k}(1+\lambda_{k})^{-\frac{\beta n}{n-1}}\right)^{n-1}
\left(\sum_{k}(1+\lambda_{k})^{-(\alpha-\beta)n}\,\E\left[ \<\,J(t) -J(s)\,,\;\phi_{k}\,\>^{2n} \right] \right) \\
&\leq&  C\,\left(\sum_{k}\frac{1}{(1+\lambda_{k})^{\frac{\beta n}{n-1}}}\right)^{n-1}
\left(\sum_{k} \frac{  (1\vee\lambda_{k})^{2n}\,\|\phi_{k}\|^{2n} }{(1+\lambda_{k})^{(\alpha-\beta)n}}\right)\,(t-s)^{n}.
\end{eqnarray*}
It follows from (\ref{E:EigenfcnUpperBound}) that (\ref{E:TightnessOfIntUNdMN}) holds true if we choose $\beta\in \left( \frac{d(n-1)}{2n},\,\alpha-\frac{d}{2}-2- \frac{d}{2n}\right)$. This choice of $\beta$ is possible if and only if $\alpha>2+d$.

\textbf{(ii) Tightness. }We will show that there exists $N_0=N_0(D)$ such that for $\alpha>d+2$,
\begin{equation}\label{E:TightnessOfIntUNdMN}
 \sup_{N>N_0} \E \left[ |J_N(t) -J_N(s)|_{-\alpha}^{4} \right] \leq  C\,(t-s)^{2}
\end{equation}
whenever $0\leq s\leq t\leq T$, where $C=C(d, D, T,\alpha, \|u_0\|, \|q\|)\,>0$ is a constant independent of $N$, $s$ and $t$. By the Kolmogorov-Centov tightness criteria (see \cite[Theorem 3.8.8]{EK86}), (\ref{E:TightnessOfIntUNdMN}) implies tightness of $\{J_N\}_{N\geq 1}$ in  $D([0,T],\,\H_{-\alpha})$.

Using H\"older inequality $(\sum_{i} a_ib_i)^n \leq (\sum_i a_i^{n/(n-1)})^{n-1}(\sum_i b_i^n)$ (with $n=2$) and the condition $\alpha>d+2$ as in step (i) above, it suffices to show that
\begin{equation}\label{E:TightnessOfIntUNdMN_phi_k}
\sup_{N\geq N_0}\E\left[ \<\,J_N(t) -J_N(s)\,,\;\phi_{k}\,\>^{4} \right]
\leq C\, (1\vee \lambda_{k})^{4}\,\|\phi_{k}\|^{4}\,(t-s)^{2}
\end{equation}
for any $0\leq s\leq t \leq T$ and $k\in \mathbb{N}$, where $N_0=N_0(D)$ and $C=C(d, D, T,\|u_0\|,\|q\|)$.

We now prove (\ref{E:TightnessOfIntUNdMN_phi_k}) by first writing
\begin{equation*}
    J_N(t) -J_N(s)= \Big(\int_0^s\mathbf{U}^N_{(t,r)}-\mathbf{U}^N_{(s,r)}\,dM^N_{r}\Big) + \int_s^t\mathbf{U}^N_{(t,r)}dM^N_{r}.
\end{equation*}
Fix $\phi_k$ and $s\leq t$. Observe that
\begin{equation*}
  \Gamma_w := \Big\langle \int_0^w\mathbf{U}^N_{(t,r)}-\mathbf{U}^N_{(s,r)}\,dM^N_{r},\;\phi_k \Big\rangle
\end{equation*}
is a martingale for $w\in [0,s]$. As in \eqref{E:JumpSizeK}, the jump size $\triangle \Gamma_w := \Gamma_w - \Gamma_{w-}$ satisfies
\begin{equation*}\label{E:JumpSizeGamma}
\sup_{w\in[0,s]}|\triangle \Gamma_w| \leq \sup_{r\in[0,s]}\|Q^N_{r,t}\phi_k- Q^N_{r,s}\phi_k\| / \sqrt{N}.
\end{equation*}
Moreover, by Theorem \ref{T:3.3}, the dual predictable projection $\< \Gamma \>$
of the quadratic variation $[\Gamma]$ of $\Gamma$ is
\begin{equation*}\label{E:QuadVarGamma}
\< \Gamma\> _w = \int^w_0 \D^{(q),N}_r(Q^N_{(r,t)}\phi_k-Q^N_{(r,s)}\phi_k)\,dr \qquad\hbox{for } w \in [0,s],
\end{equation*}
where $\D^{(q),N}_r(\phi,\psi) := \< \textbf{a}\nabla\phi\cdot\nabla\psi+q_N(r)\phi\psi,\, \X^N_r\>$ and $\D^{(q),N}_r(\phi) := \D^{(q),N}_r(\phi,\phi)$.
Therefore, by Burkholder-Davis-Gundy inequality for discontinuous martingales (see the remark after Theorem 74 in Chapter IV  of \cite{peP05}), we have $\E[\,\Gamma_s^4\,]\leq \bar{c}\,\E[\,[\Gamma]_{s}^2\,]$ for some absolute constant $\bar{c}$. Hence, argue as in \eqref{e:4.30}, and then by Lemma \ref{L:BoundVariance_stphi} and Lemma \ref{L:Bound_QtMinusQs}, we obtain
\begin{eqnarray*}
&& \E[\,\Gamma_s^4\,] \leq \bar{c}\,\E[\,[\Gamma]_{s}^2\,]\\
&\leq&  2\bar{c} \,\bigg(\,\E[\,\<\Gamma\>_s^2\,] + \frac{\sup_{r\in[0,s]}\|Q^N_{r,t}\phi_k- Q^N_{r,s}\phi_k\|^2 }{N}\,\E[\<\Gamma\>_s] \,\bigg)\\
&\leq& 2\bar{c}\,\E[\,\<\Gamma\>_s^2\,] \,+\,C\,\frac{\|\phi_k\|^2\Big(\lambda_k^2(t-s)^2+(t-s)\Big)}{N}\Big(\lambda_k^2+\|\phi_k\|^2+\lambda_k^2\|\phi_k\|^2\Big)(t-s).
\end{eqnarray*}
where $C=C(D,T,\|u_0\|,\|q\|)>0$
Estimating $\E[\,\<\Gamma\>_s^2\,]$ by the argument we used for $H^k_N(t)$ in the proof of Lemma \ref{L:HNGNtoZero} (via an inequality in Remark \ref{Rk:BoundVariance_stphi}), we see that
$$\E\left[ \Big\langle \int_0^s\mathbf{U}^N_{(t,r)}-\mathbf{U}^N_{(s,r)}\,dM^N_{r},\;\phi_k \Big\rangle^{4} \right]=\E[(\Gamma_s)^4]$$
is bounded above by the RHS of (\ref{E:TightnessOfIntUNdMN_phi_k}) for $N\geq N_0(D)$.

Similarly, by consider the martingale
\begin{equation*}
  \Theta_w := \Big\langle \int_s^{s+w}\mathbf{U}^N_{(t,r)}dM^N_{r},\;\phi_k \Big\rangle, \quad w\in[0,t-s];
\end{equation*}
and by using Lemma \ref{L:BoundVariance_stphi}, we can check that $\E\Big[ \< \int_s^{t}\mathbf{U}^N_{(t,r)}\,dM^N_{r},\;\phi_k\>^{4} \Big]= \E[(\Theta_{t-s})^4]$ is bounded above by the RHS of (\ref{E:TightnessOfIntUNdMN_phi_k}) for $N\geq N_0(D)$.

\textbf{(iii) Convergence of finite dimensional distributions. }
As in the proof of Theorem \ref{T:3.8}, it suffices to show that  as $N \rightarrow \infty$,
\begin{equation}\label{E:FiniteDimConvergence_Integral}
(\<J_N(t_1),\,\psi_1\>, \dots, \<J_N(t_n),\,\psi_n\>) \, \toL \,(\<J(t_1),\,\psi_1\>,\dots,\<J(t_n),\,\psi_n\>) \quad\text{in } \R^n
\end{equation}
for any $n\in \mathbb{N}$, $0\leq t_1\leq\cdots \leq t_n<\infty$ and $\{\psi_j\}_{j=1}^n\subset C(\bar{D})$.

For $n=1$, fix $t\geq 0$ and $\phi\in C(\bar{D})$. Note that  $\theta \mapsto \<\int^\theta_0\mathbf{U}^N_{(t,s)}dM^N_{s} ,\, \phi\>$ is a martingale for $\theta \in [0,t]$, with jumps size at most
$\sup_{\theta\in[0,t]}|M^N_\theta(Q^N_{(\theta,t)}\phi)-M^N_{\theta-}(Q^N_{(\theta,t)}\phi)| \leq \|\phi\|/\sqrt{N}$,
by (\ref{E:ContractionQNQ}) and (\ref{E:JumpSizeM^N}). Hence by the functional central limit theorem for real-valued martingales (see \cite{LS80}),
\begin{equation}
\bigg\{ \big\langle \int^\theta_0 \mathbf{U}^N_{(t,s)}dM^N_{s} ,\, \phi \big\rangle ; \ \theta\in[0,t] \bigg\} \,\toL \, \bigg\{ \big\langle \int^\theta_0\mathbf{U}_{(t,s)}dM_{s},\,\phi\big\rangle ; \ \theta\in[0,t] \bigg\} \quad\text{in } D([0,t],\R)
\end{equation}
as $N \rightarrow \infty$.

For an integer $n>1$, (\ref{E:FiniteDimConvergence_Integral}) follows from Lemma \ref{L:ConvergenceDistOfIntegral} and the towering property
\begin{equation*}
\E Z=\E \E[Z|\F_{t_1}]= \E \E[\E[Z|\F_{t_2}]|\F_{t_1}]= \cdots \quad\text{for }0\leq t_1\leq t_2 \leq t_3 \leq \cdots.
\end{equation*}
We illustrate this for the case $n=3$; the proof for the general case is the same.
The Fourier transform
\begin{eqnarray*}
 \E\bigg[\exp{\Big(i \sum^3_{k=1}a_k \<J^N_{t_k},\,\psi_k\>\Big)}\bigg]
&=&
\E\Bigg[
    \exp{
        \Big(
            i \sum^3_{j=1} a_j \Big\langle
                \int^{t_1}_0 \mathbf{U}^N_{({t_j},s)}dM^N_{s},\,\psi_j
            \Big\rangle
        \Big)
    } \,
   \\  &&\quad\cdot
    \E\bigg[
        \exp{
            \Big(
                i \sum^3_{j=2} a_j \Big\langle
                    \int^{t_2}_{t_1} \mathbf{U}^N_{({t_j},s)}dM^N_{s},\,\psi_j
                \Big\rangle
            \Big)
        }
       \\ &&\quad\quad \cdot
        \E\Big[
            \exp{
                \Big(
                    i a_3 \Big\langle
                        \int^{t_3}_{t_2} \mathbf{U}^N_{({t_3},s)}dM^N_{s},\,\phi_3
                    \Big\rangle
                \Big)
            }
        \Big|\F_{t_2}\Big]
    \bigg|\F_{t_1}\bigg]
\Bigg].
\end{eqnarray*}
We then apply Lemma \ref{L:ConvergenceDistOfIntegral} three times successively, starting from the inner most term involving $\F_{t_2}$. Hence we have convergence (\ref{E:FiniteDimConvergence_Integral}).

The proof of the lemma is complete.
\end{pf}

\subsection{Characterization of $\Y$}

Let $\Y$ be any subsequential limit of $\Y^N$.
By Theorem \ref{T:3.7}, $\Y$ has a continuous version in $\H_{-\alpha}$
for every $\alpha > d \vee (d/2+2)$.
It follows from Theorems \ref{T:3.3}, \ref{T:3.8} and \ref{T:ConvergenceStocInt} that we have
\begin{equation}\label{e:3.45}
\Y_t \,  \equalL \, \mathbf{U}_{(t,0)}\Y_0 + \int_0^t \mathbf{U}_{(t,s)}\,dM_s,\quad \text{in }D([0,\,T],\H_{-\alpha}).
\end{equation}

\begin{thm}\label{T:Convergence_RBM_Y}
   The limiting process $\Y$ is a continuous Gaussian Markov process that is unique in distribution. Moreover, $\Y$ has a version in $C^{\gamma}([0,\infty),\,\H_{-\alpha})$ for $\gamma \in (0,1/2)$.
\end{thm}

\begin{pf}
Since $M$ is Gaussian, $\int_0^t \mathbf{U}_{(t,s)}\,dM_s$ is a Gaussian process by the construction of the stochastic integral. On the other hand, $\mathbf{U}_{(t,0)}\Y_0$ is a Gaussian process and is independent of $\int_0^t \mathbf{U}_{(t,s)}\,dM_s$ since $M$ has independent increments. Therefore $\Y_t$, as the sum of two independent Gaussian processes, is a Gaussian process.

The Markov property of $\Y$ is basically due to the independent increments of the differentials; see section 5.6 of \cite{peP05}.
For reader's convenience, we give a proof that $\Y$ is a Markov process with respect to  its own filtration
$\F^{\Y}_t:= \sigma(\Y_r:\,r\leq t)= \sigma(\<\Y_r,\phi\>:\,r\leq t,\,\phi\in\H_{\alpha})$.
We in particular have from \eqref{e:3.45} that for $s\leq t$,
    \begin{equation}\label{E:GeneralizedOUformula_RBM_Killing_st}
        \Y_t \equalL \mathbf{U}_{(t,s)}\Y_s + \int_s^t \mathbf{U}_{(t,r)}\,dM_r \quad \text{ in } \H_{-\alpha}.
    \end{equation}
Together with the fact that $M$ has independent increments, we have
    \begin{equation}\label{E:CovraianceStructure_U}
       {\rm Cov}(\<\Y_s,\phi\>,\, \<\Y_t,\psi\>)= {\rm Cov} (\<\Y_s,\phi\>,\, \<\mathbf{U}_{(t,s)}\Y_s,\psi\>)
    \end{equation}
for all $s\leq t$ and $\phi,\,\psi\in \H_{\alpha}$. To show that $\Y$ is Markov, note that (\ref{E:CovraianceStructure_U}) together with the fact that $\mathbf{U}_{(t,s)}\Y_s\in \F^{\Y}_s$ yield
$\E[F(\Y_t)\big|\F^{\Y}_s]= F\big(\mathbf{U}_{(t,s)}\Y_s\big)$ for all $F\in C_b(\H_{-\alpha})$. Using (\ref{E:GeneralizedOUformula_RBM_Killing_st}) and the fact that $\mathbf{U}_{(t,s)}\Y_s\in \sigma(\Y_s)$, we obtain $\E[F(\Y_t)\big|\Y_s]= F\big(\mathbf{U}_{(t,s)}\Y_s\big)$ for all $F\in C_b(\H_{-\alpha})$. This shows that  $\Y$ is Markov.

The H\"older continuity of $\Y$ follows immediately from Theorem \ref{T:3.8} and Theorem \ref{T:ConvergenceStocInt}.
\end{pf}

\medskip

The proof of Theorem \ref{T:Convergence_RBM_killing_Y} is now complete.

\section{Appendix}

\subsection{Hilbert-Schmidt Operators}

 Hilbert-Schmidt  operators appear  naturally in stochastic analysis in infinite dimensions. The main properties of these operators can be found in standard references (e.g. \cite{GT95}). We now recall the main definitions.

\begin{definition}
Let $X=(X_t)_{t\geq 0}$ be an $\H_{-\alpha}$-valued process defined on a probability space $(\Omega,\,\mathcal{F},\,\P)$. We say $X$ is (centered) \textbf{Gaussian} if $\{X_t(\phi):\,\phi\in\H_{\alpha},\,t\in[0,\infty)\}$ form a (centered) Gaussian system. That is, $\left( X_{t_1}(\phi_1),\cdots,\,X_{t_k}(\phi_k) \right)$ is a (centered) Gaussian vector in $\R^k$ for any $k\in \mathbb{N}$, any $\{t_i\}_{i=1}^k \subset [0,\infty)$ and any $\{\phi_i\}_{i=1}^k\subset \H_{\alpha}$. We say $X$ is \textbf{continuous} if $t\mapsto X_t$ is continuous $\P$-a.s. $X$ is said to be \textbf{square-integrable} if $\E[|X_t|^2_{-\alpha}]<\infty$ for all $t\geq 0$. Finally, we say $\X$ has \textbf{independent increments} if for any $0\leq s<t$ and $\phi\in\H_{\alpha}$, the real random variable $X_t(\phi)-X_s(\phi)$ is independent of the $\sigma$-field generated by $\{X_r(\psi):\,0\leq r\leq s,\,\psi\in\H_{\alpha}\}$.
\end{definition}

Suppose $X$ and $Y$ are real separable Hilbert spaces with inner product $\<\,,\,\>_X$ and $\<\,,\,\>_Y$ (we simply write $\<\,,\,\>$ when there is no confusion for which Hilbert space we are considering). The class of bounded linear operators from $X$ to $Y$ will be denoted by $L(X,Y)$ ($L(X)$ for short when $X=Y$). It is well known that $A\in L(X,Y)$ is \textbf{compact} (i.e. the range of the unit sphere in $X$ is relatively compact in $Y$) if and only if there exist orthonormal systems (ONS for short) $\{e_n\}\subset X$, $\{f_n\}\subset Y$ and a sequence of real numbers $a_n\to 0$ such that $A$ has the representation
\begin{equation}\label{E:RepresentationCompactOperator}
Ax=\sum_{n\geq 1} a_n\,\<x,\,e_n\>\,f_n \quad \text{for all }x\in X.
\end{equation}

\begin{definition}
\begin{enumerate}
\item   $A\in L(X,Y)$ is said to be \textbf{Hilbert-Schmidt} (denoted by $A\in L_2(X,Y)$) if $A$ has the representation  \eqref{E:RepresentationCompactOperator} with $\sum_{n\geq 1} a_n^2<\infty$. In this case, the \textbf{Hilbert-Schmidt norm} of $A$ is defined to be
    $$\|A\|_2:= \left(\sum_{n\geq 1} a_n^2\right)^{1/2} = \left(\sum_{n\geq 1} |Ae_n|^2 \right)^{1/2} $$
    Note that $\|A\|_2$ is independent of the choice of the ONS $\{e_n\}\subset X$.

\item   The \textbf{Trace} of $A \in L(X)$ is
    $$Tr(A):= \sum_{n\geq 1} \<Ae_n,\,e_n\>$$
    Note that $Tr(A)$ is independent of the choice of the ONS $\{e_n\}\subset X$.
\end{enumerate}
\end{definition}

The following lemma is equivalent to the statement that $(\Phi_{imb},\,\H_{\beta},\,\H_{\gamma})$ is an abstract Wiener space if $\beta>d/2+\gamma$ (cf. \cite{sS07}).

\begin{lem}
For any $\beta,\,\gamma\in\R$ with $\beta>\gamma +d/2$, the imbedding $\Phi_{imb} :\;\H_{\beta}\rightarrow \H_{\gamma}$ is Hilbert-Schmidt.
\end{lem}
\begin{pf}
We want to show that $\sum_{k}\Big|\Phi_{imb}\left(h^{(\beta)}_{k}\right)\Big|_{\gamma}^2 <\infty$. The left hand side equals
$$\sum_{k} (1+\mu_{k})^{-\beta}\big|\phi_{k}\big|_{\gamma}^2= \sum_{k} (1+\lambda_{k})^{-\beta+\gamma}.$$ By Weyl's formula (\ref{E:WeylLaw}), the latter quantity is finite if and only if $$\int_{1}^{\infty} (1+x)^{-\beta+\gamma}x^{d/2-1}\,dx <\infty.$$
This is true if and only if $\beta-\gamma>d/2$.
\end{pf}

\subsection{$\int_0^{t} \mathbf{U}_{(t,s)}dM_{s}$ is well defined}\label{subsection:WellDefineStocInt}

As mentioned earlier, we have to make sure that $\mathbf{U}_{(t,s)}$ (for $s\in[0,t]$) lies within the class of integrands with respect to  $M$. We will follow the construction of stochastic integrals with respect to  Hilbert space valued r.c.l.l. square-integrable martingales in \cite{MP80}. See \cite{pDjZ08, GT95, peP05} for more comprehensive and recent treatments.

We denote by $M_c^2([0,\infty),\H_{-\alpha})$  the class of continuous square-integrable $\H_{-\alpha}$-valued martingales with zero initial value. Fix $\alpha>d\vee (d/2+1)$ and recall  from Theorem \ref{L:ConvergenceOfM^N_RBM_killing}  that $M\in M_c^2([0,\infty),\H_{-\alpha})$ is Gaussian, has independent increments and  covariance
    \begin{equation}
    \tilde{\E} \left[ \<M_s,\phi\>\<M_t,\psi\> \right] = \int_0^{s\wedge t} \D^{(q)}_r(\phi,\psi)\,dr,
    \end{equation}
where $\D^{(q)}_r$ is the bilinear form on $\H_{-\alpha}$ defined in (\ref{Def:D_r(phipsi)}). We will omit the filtration when there is no ambiguity. For example, we simply say that $M$ is adapted rather than $\tilde{\mathcal{F}}_{t}$-adapted since it is defined on $(\tilde{\Omega},\,\tilde{\mathcal{F}},\,\tilde{\mathcal{F}}_{t},\,\tilde{\P})$.
For $T\in (0, \infty]$, denote by $\mathfrak{P}_{[0, T]}$   the $\sigma$-field of predictable sets on $\tilde{\Omega}\times [0,T]$. That is, the smallest $\sigma$-field  making all adapted processes with left continuous paths measurable (c.f. p.156 of \cite{peP05} or section 1.7 \cite{MP80}).
When $T=\infty$, we write $ \mathfrak{P} $ for $\mathfrak{P}_{[0, \infty)}$.

By a direct calculation,
\begin{equation}
[[M]]_t:= \sum_{k}\int_0^t \D^{(q)}_r(h_{k}^{(\alpha)})\,dr
\end{equation}
is the unique continuous, adapted and increasing real process such that $\big|M_t\big|^2_{-\alpha}-[[M]]_t$ is a real martingale (cf. Remark 2.2 in \cite{GT95}). $[[M]]$ is called the \textbf{real increasing process} associated to $M$. Besides, the operators $Q_s:\,H_{-\alpha}\rightarrow H_{-\alpha}$ (for $s\geq 0$) defined by
\begin{equation}
\big\langle Q_s(h^{(-\alpha)}_{i}) ,\,h^{(-\alpha)}_{j}\big\rangle_{-\alpha} := \dfrac{\D^{(q)}_s(h^{(\alpha)}_{i},\,h^{(\alpha)}_{j})}{\sum_{k}\D^{(q)}_s(h^{(\alpha)}_{k})}
\end{equation}
is called the \textbf{characteristic operator process} associated to $M$. Clearly, $Q_s$ is a non-negative operator on $H_{-\alpha}$ with $Tr(A)=1$ where `Tr' means `Trace'. As a remark, the operator-valued process $\<\<M\>\>_t := \int_0^t Q_s\,d[[M]]s$ (in the sense of Bochner's integral) is called the \textbf{operator increasing process associated to} $M$ and plays an analogous role as the  quadratic variation  of real-valued martingales (see Theorem 2.3 in Chapter 1 of \cite{GT95} for its basic properties).

Following \cite{MP80}, the class of possible integrands for the stochastic integral with $M_t$ as integrator (on the interval $[0,T]$) can be defined as follows: On the space of $\mathfrak{P}_{[0, T]}$-simple $L(\H_{-\alpha})$-valued processes, we define a scalar product
\begin{equation}\label{E:ScalarProductSImpleProcesses}
    (A,B):= \tilde{E}\left[\int_0^T Tr(A\,Q_s\,B^{\ast})\,d[[M]]_s\right],
\end{equation}
where $B^{\ast}$ is the adjoint of the operator $B$. The completion of the $\mathfrak{P}_{[0, T]}$-simple $L(\H_{-\alpha})$-valued processes with respect to  the scalar product in (\ref{E:ScalarProductSImpleProcesses}), denoted by $\Lambda^2(\H_{-\alpha},\mathfrak{P}_{[0, T]},M)$, is the desired class of integrands. It is worth noting that (c.f. p.171 \cite{MP80}) $\Lambda^2(\H_{-\alpha},\mathfrak{P}_{[0, T]},M)$ contains processes whose values may be unbounded operators.

By section 1.3 of \cite{GT95}, $\Lambda^2(\H_{-\alpha},\mathfrak{P}_{[0, T]},M)$ contains the class of all processes $(\Phi_t)_{t\in[0,T]}$ such that
\begin{enumerate}
\item[(i)]  $\Phi_t$ is a linear operator (not necessarily bounded) from $\sqrt{Q_t}\,\H_{-\alpha}$ to $\H_{-\alpha}$ such that $\Phi_t\sqrt{Q_t}\in L_2(\H_{-\alpha})$ is Hilbert-Schmidt for all $t\in[0,T]$ a.s.,
\item[(ii)] $\Phi_t\sqrt{Q_t}$ is $\mathfrak{P}\big|_{\tilde{\Omega}\times [0,T]}$-measurable (i.e. predictable), and
\item[(iii)] $\E\left[\int_0^T \|\Phi_t\sqrt{Q_t}\|_2^2\,d[[M]]_t \right]<\infty$
            where $\|\cdot\|_2$ is the Hilbert-Schmidt norm.
\end{enumerate}

Now for any $t>0$, the deterministic process $(U_{(t,\theta)})_{\theta\in[0,t]}$ lies in the class of integrands with respect to  $M$. This is because on one hand
\begin{eqnarray*}
\|U_{(t,\theta)}\sqrt{Q_{\theta}}\|_2^2 &=& Tr\left(U_{(t,\theta)}Q_{\theta}U_{(t,\theta)}^{\ast}\right) \quad\text{ the trace of }U_{(t,\theta)}Q_{\theta}U_{(t,\theta)}^{\ast}\\
&=& \sum_{k}\<U_{(t,\theta)}Q_{\theta}U_{(t,\theta)}^{\ast}(h^{(-\alpha)}_{k}),\,h^{(-\alpha)}_{k}\>_{-\alpha}\\
&=& \sum_{k}\<Q_{\theta}U_{(t,\theta)}^{\ast}(h^{(-\alpha)}_{k}),\,U_{(t,\theta)}^{\ast}(h^{(-\alpha)}_{k})\>_{-\alpha}\\
&=& \dfrac{\sum_{k}\D^{(q)}_{\theta}(Q_{(\theta,t)}h^{(\alpha)}_{k})}{\sum_{i}\D^{(q)}_{\theta}(h^{(\alpha)}_{i})}
\end{eqnarray*}
which is finite provided that $u_0$ is not identically zero; and on the other hand, by Lemma \ref{L:BoundVariance_stphi},
\begin{eqnarray*}
&& \E\left[\int_0^t \|U_{(t,\theta)}\sqrt{Q_{\theta}}\|_2^2\,d[[M]]_{\theta} \right]\\
&=& \sum_{k} \int_0^t\D^{(q)}_{\theta}(Q_{(\theta,t)}h^{(\alpha)}_{k})\,d\theta \\
&\leq& C(d,D,T)\,\|u_0\|\,\left(t\,\sum_{k}\frac{\lambda_{k}+\|\phi_{k}^2\|}{(1+\lambda_{k})^{\alpha}} + \|q\|\,t^{3/2}\,\sum_{k}\frac{\|\phi_{k}\|^2}{(1+\lambda_{k})^{\alpha}} \right)\quad\text{for }t\in[0,T]\\
&<&\infty \quad \text{if  }\alpha>d \vee (d/2+1).
\end{eqnarray*}

We conclude that for any fixed $t\geq 0$, $\left\{\int_0^{s}U_{(t,\theta)}\,dM_{\theta};  s\in[0,t]\right\}$ is a continuous, adapted square-integrable $\H_{-\alpha}$-valued martingale with  $\E\left[\big|\int_0^{s}U_{(t,\theta)}\,dM_{\theta}\big|_{-\alpha}^2\right]=\sum_{k} \int_0^s\D^{(q)}_{\theta}(Q_{(\theta,t)}h^{(\alpha)}_{k})\,d\theta$. In particular, putting $s=t$, we have that $\int_0^{t}U_{(t,\theta)}\,dM_{\theta}$ is a well defined $\tilde{\mathcal{F}}_{t}$-measurable $\H_{-\alpha}$-valued random variable with finite second moment. Moreover, since $M$ is centered Gaussian with independent increments,   $\int_0^{t}U_{(t,\theta)}\,dM_{\theta}$ is also centered Gaussian.

\subsection{An identity}

The following equality is used in Lemma \ref{L:Bound_QtMinusQs}.
\begin{lem}\label{L:Gamma_Half_k}
    \begin{equation*}
        \int_0^s\cdots \int_0^{s_{k}}
        \dfrac{1}{\sqrt{(s-s_2)(s_2-s_3)\cdots (s_k-s_{k+1})}}
        \,ds_{k+1} \cdots \,ds_2
        = \dfrac{\pi^{k/2}}{\Gamma\left(\frac{k+2}{2}\right)}\,s^{k/2}.
    \end{equation*}
\end{lem}

\begin{pf} Denote the integral on the left hand side as $V_k$.
    For any $a\in(0,\infty)$,
    $$\int_0^x \frac{y^a}{\sqrt{x-y}}\,dy= \dfrac{\sqrt{\pi}\,\Gamma(1+a)}{\Gamma(3/2+a)}\,x^{1/2+a}$$
     Using this, we can iterate it to obtain $V_k=\int_0^t c_k\,s^{k/2}\,ds$,
      where
    $$c_1=2 \quad\text{and}\quad c_{k+1}=c_k\,\dfrac{\sqrt{\pi}\,\Gamma(1+k/2)}{\Gamma(3/2+k/2)} \text{ for }k\geq 2.$$
\end{pf}

\subsection{Reflected diffusions killed by local time}\label{S:5.4}

Suppose now, instead of being killed by $q_N$, that $Z^{(N)}_i=Z_i$ is the subprocess of $X_i$ killed by $2\int_0^tq(s,X_i(s))dL^{(i)}_s$ \emph{for all} $N$. In Remark \ref{R:1.7}(ii), we claimed that Theorem \ref{T:HydrodynamicLimit_RBM_Killing} and Theorem \ref{T:Convergence_RBM_killing_Y} remain valid. The claim that Theorem \ref{T:HydrodynamicLimit_RBM_Killing} remains true is easy to be verified. We now provide some details to support the claim that Theorem \ref{T:Convergence_RBM_killing_Y} remains valid.

By the same proof of Lemma \ref{L:KeyMtgRobinModel}, we have the following:
\begin{lem}\label{L:KeyMtgRobinModel_LocalTime}
Fix any positive integer $N$. For any $\phi\in Dom^{Feller}(\A)$, we have under $\P^{\mu}$ for any $\mu\in E_N$,
\begin{eqnarray*}
M^{\phi}_t &:=& \<\phi,\X^N_t\>-\<\phi,\X^N_0\>-\int_0^t \<\A \phi,\,\X^N_s\>\,ds + \frac{1}{N}\sum_{i=1}^N\int_0^tq(s,Z_i(s))\phi(Z_i(s))\,dL^i_s \label{E:MtgMphi_LocalTime}
\end{eqnarray*}
is an $\F^{\X^N}_t$-martingale under $\P^{\mu}$ for any $\mu\in E_N$. Moreover, $M^{\phi}_t$ has predictable quadratic variation
\begin{equation*}\label{E:QuadVar_MtgMphi_LocalTime}
\<M^{\phi}\>_t = \frac{1}{N}\,\bigg[\,\int_0^t\< \textbf{a}\nabla \phi\cdot  \nabla \phi,\,\X^N_s\>\,ds \,+\,
\frac{1}{N}\sum_{i=1}^N\int_0^tq(s,Z_i(s))\phi^2(Z_i(s))\,dL^i_s \,\bigg].
\end{equation*}
Moreover, (\ref{E:QuadVar_MtgMphi_2}) still holds for this new martingale.
\end{lem}

Starting from the above lemma, we just need slight modifications in the proof of Theorem \ref{T:Convergence_RBM_killing_Y}. It is easier in this case since now we have  $Q^N=Q$ and $\mathbf{U}^N=\mathbf{U}$. Note that in the proof of Lemma \ref{L:ConvergenceDistOfIntegral}, the expressions $H_N(t)$ and $G_N(t)$ in (\ref{E:ForthTerm}) should be replaced by, respectively,
$$ \frac{1}{N}\sum_{i=1}^N\int_0^t\psi^2(Z^i_r)q(r,Z^i_r)\,dL^i_r \quad\text{and}\quad \int_0^t\int_{\partial D}\psi^2(z)q(r,z)\,u(r,z) \rho(z)\,d\sigma(z)\,dr.$$
In addition, we  should also use the following lemma rather than Lemma \ref{L:HNGNtoZero}.

\begin{lem}\label{L:HNGNtoZero_LocalTime}
Let $\{\phi(r)\,:\,r\geq 0\}\subset C(\bar{D})$ be such that $\sup_{r\geq 0}\|\phi(r)\| <\infty$. For any $p\geq 1$, we have
\begin{equation*}
\lim_{N\to\infty} \E\left[\,\bigg(\sup_{t\in[0,T]}\Big|\sum_{i=1}^N\int_0^t \phi(r,Z^i_r)\,dL^i_r\,-\,
\int_0^t\int_{\partial D}\phi(r,z)\,u(r,z)\,\rho(z)\,d\sigma(z)\,dr\Big|\bigg)^p\,\right] =0.
\end{equation*}
\end{lem}

The proof of Lemma \ref{L:HNGNtoZero_LocalTime}  is the same as that of Lemma \ref{L:HNGNtoZero}.

\vspace{5mm}
    \textbf{Zhen-Qing Chen}

    Department of Mathematics, University of Washington, Seattle, WA 98195, USA

    Email: zqchen@uw.edu
\vspace{2mm}

    \textbf{Wai-Tong (Louis) Fan}

    Department of Mathematics, University of Washington, Seattle, WA 98195, USA

    Email: louisfan@math.washington.edu

\end{document}